\documentclass[11pt,pedantic]{volume}
\usepackage{graphicx}

\usepackage{latexsym}
\usepackage{amsfonts}
\usepackage{amsthm}
\usepackage{amsmath}

\newcommand{\reals}{\mathbb{R}}
\let\eps\varepsilon

\newtheorem{theorem}{Theorem}[section]
\newtheorem{lemma}[theorem]{Lemma}
\newtheorem{proposition}[theorem]{Proposition}
\newtheorem{corollary}[theorem]{Corollary}
\newtheorem*{theorem*}{Theorem}

\newcommand{\lemlab}[1]{\label{lemma:#1}}
\newcommand{\theolab}[1]{\label{theo:#1}}
\newcommand{\figlab}[1]{\label{fig:#1}}
\newcommand{\seclab}[1]{\label{section:#1}}
\newcommand{\lemref}[1]{\ref{lemma:#1}}
\newcommand{\theoref}[1]{\ref{theo:#1}}
\newcommand{\figref}[1]{\ref{fig:#1}}
\newcommand{\secref}[1]{\ref{section:#1}}

\def\complaint#1{}
\def\withcomplaints{
\newcounter{mycomplaints}
\def\complaint##1{\refstepcounter{mycomplaints}%
\ifhmode%
\unskip%
{\dimen1=\baselineskip \divide\dimen1 by 2 %
\raise\dimen1\llap{\tiny -\themycomplaints-}}\fi%
\marginpar{\tiny [\themycomplaints]: ##1}}%
}

\newcommand{\phce}{\bar X}
\newcommand{\xcone}{\bar X_0}
\newcommand{\pce}{X}
\newcommand{\eeqref}[1]{\eqref{eq:#1}}
\newcommand{\im}{\operatorname{Im}}
\newcommand{\Area}{\operatorname{Area}}

\withcomplaints

\begin{document}

%% Title, optional running title in brackets.
\title{{%
Expansive Motions and the \\Polytope of Pointed Pseudo-Triangulations}}

\runningtitle
{Expansive Motions and the Pseudo-Triangulation Polytope}

%% Later we might need a title for the Table of Contents of the entire
%%  volume. So put this in please.
\toctitle{Expansive motions and the polytope of pointed
  pseudo-triangulations}

%% Stuff that usually goes into \thanks is moved toward the end of the
%% paper.  Do NOT use \thanks.
% \thanks{%
%        Work on this paper has been supported .}

%% Authors do not have affiliations next to them, affiliations go to
%  the end of the paper.
\author{
G\"unter Rote \and
Francisco Santos \and
Ileana Streinu}

%% similar to runningtitle/toctitle
\tocauthor
{G\"unter Rote, \and
Francisco Santos, \and
Ileana Streinu}
\runningauthor{G.~Rote, \and
F.~Santos, \and
I.~Streinu}

\maketitle

\begingroup
\let\eps\varepsilon
\mathcode`w=\omega

\begin {abstract}

We introduce the polytope of pointed pseudo-triangulations of a point set in
the plane, defined as the polytope of infinitesimal expansive motions
of the points subject to certain constraints
on the increase of their distances. Its $1$-skeleton is the graph
whose vertices are the pointed
pseudo-triangulations of the point set and whose edges are flips of
interior pseudo-triangulation edges.

For points in convex position we obtain a new realization
of the associahedron, i.e., a geometric representation of
the set of triangulations of an $n$-gon, or of the set of binary
trees on $n$ vertices, or of many other combinatorial objects
that are counted by the Catalan numbers.
By considering the $1$-dimensional version of the polytope of constrained
expansive motions we obtain a second  distinct realization of the
associahedron as a perturbation of the positive cell in a Coxeter
arrangement. 

Our methods produce as a by-product a new proof that every simple polygon or
polygonal arc in the plane has expansive motions, a key step in the proofs of
the Carpenter's Rule Theorem by
Connelly, Demaine and Rote (2000) and by Streinu~(2000).
\end{abstract}

\section{Introduction}
\label{introduction}

\paragraph{\bf Polytopes for combinatorial objects.}
Describing all instances of a combinatorial structure (e.g.\ trees or
triangulations) as vertices of a polytope is often a step towards
giving efficient optimization algorithms on those structures. It
also leads to quick prototypes of enumeration algorithms
using known vertex enumeration techniques
and existing code \cite{avis_fukuda,gdd}.

One particularly nice example is the associahedron,
(see Figure~\ref{fig:assoc} for an example):
the vertices of this polytope correspond to Catalan structures.
The Catalan structures refer to any of a great number of
combinatorial objects which are counted by the Catalan numbers
(see the extensive list in Stanley~\cite[ex.~6.19, p.~219]{Stanley2}).
Some of the most notable  ones are
the triangulations of a convex polygon, binary trees,
the ways of evaluating a product of $n$ factors when multiplication
is not associative (hence the name associahedron), and
monotone lattice paths that go from one corner of a square to the
opposite corner without crossing the diagonal.

In this paper we describe
a new polyhedron whose vertices correspond to 
\emph{pointed pseudo-triangulations}.

\paragraph {\bf Pseudo-triangulations.}
Pseudo-triangulations, as well as  the closely related geo\-desic triangulations
of simple polygons,
have been used in Computational Geometry in
applications such as visibility \cite{PV,PV2,PV-theory,Toth-Speck}, 
ray shooting \cite{gt}, and kinetic data structures
\cite{guibas,KSS}.
The \emph{minimum} or \emph{pointed
pseudo-triangulations} introduced in Streinu~\cite{Streinu-2000}
have applications
to non-colliding motion planning of planar robot arms.
They also have very nice combinatorial and rigidity theoretic properties. The
polytope we define in this paper adds to the former, and is constructed
exploiting the latter.

\paragraph {\bf Expansive motions.}
An expansive motion on a set of points $P$
is an infinitesimal motion of the points
such that no distance between them decreases.

Expansive motions were instrumental in the first proof of
the Carpenter's Rule Theorem by
Connelly, Demaine and Rote~\cite{CDR}:
Every simple polygon or polygonal arc in the plane can be
unfolded into convex position without collisions.
Streinu~\cite{Streinu-2000} built on this work,
realizing the importance of pseudo-triangulations
in connection with expansive motions
and studying their rigidity properties.
This paper provides a systematic study of expansive motions in 
one and two dimensions.
The expansive motions of a set of $n$ points in
the plane form a polyhedral cone of dimension $2n-3$
(the \emph{expansion cone}).
As by-products of our approach we get a new proof of the existence
of expansive motions for non-convex polygons and polygonal arcs
(Theorem \theoref{unfold-ex}) and a characterization of the extreme rays of
the expansion cone of a planar point set in general position, as
equivalence classes of pointed pseudo-triangulations with one convex hull
edge removed, modulo rigid subcomponents (Proposition~\lemref{ex-rays}).

Our tool is the introduction of
\emph{constrained expansions} as expansive motions
with a special lower bound on the edge length increase.
They form a polyhedron obtained by translation of the facets of
the expansion cone. Our main result is the following (see a more precise
statement as Theorem~\ref{main}):

\begin{theorem*}
Let $P$ be a set of $n$ points in general position in the plane, $b$ of them
in the boundary of the convex hull. Then,
there is a choice of constraints 
which produces as constrained expansions of $P$ a simple polyhedron of
dimension $2n-3$ with a unique maximal bounded face of dimension 
$2n-b-3$ whose vertices correspond
to pointed pseudo-triangulations and edges correspond to flips between
them.
\end{theorem*}

The flips mentioned in the statement are a certain neighborhood
structure among pointed pseudo-triangulations
(flips of interior edges). See the details in Section~\ref{preliminaries}. 

\paragraph {\bf Two appearances of the associahedron.}
For points in convex position,
pseudo-triangulations coincide with triangulations.
We prove (Corollary \ref{coro:convex})
that, in this case, our construction gives a polytope affinely equivalent
to the standard $(n-3)$-dimensional associahedron
obtained as a secondary polytope of the point
set~\cite[Section~9.2]{Ziegler}. Perhaps surprisingly, this shows that the
secondary polytope of $n$ points in convex position in the plane
(which lives in $\reals^n$) can naturally be
embedded as a face in a $(2n-3)$-dimensional unbounded polyhedron.

The associahedron appears again as the analog of our construction for points
in one dimension (Section~\ref{1dim}).

\paragraph{\bf Rigidity.}
The connection of these results with rigidity theory is also worth mentioning.
Pointed pseudo-triangulations are special instances of infinitesimally
minimally rigid frameworks in dimension $2$, whose combinatorial
structure is well understood (see \cite{gss}).
One-dimensional minimally rigid frameworks are trees, another
well understood combinatorial structure. Adding the
constraint of expansiveness is what leads to pointed pseudo-triangulations
in $2$d, and to the special non-crossing and alternating
trees which appear in Section~\ref{1dim}.

\paragraph{\bf Future perspectives.}
It is our hope that the insight into one- and two-dimensional
motions may eventually lead to generalizations
to higher dimensions.
There is no satisfactory definition of an analog of pseudo-triangulations in
$3$ dimensions.
The $3$-dimensional version of the robot arm motion planning
problem, with potential applications to computational biology 
(protein folding), is much more challenging.

\paragraph{\bf Overview.}
In Section \ref{preliminaries}
we give the preliminary definitions and results.
Section \ref{polytope} contains the main
result, the construction of the polytope of pointed pseudo-triangulations
(\emph{ppt-polytope}).
Section~\secref{appl} applies the main result to get a
new proof for the existence of expansive motions for non-convex
polygons and polygonal arcs in the plane.
In Section \ref{other}
we present an alternative construction
of the ppt-polytope and two special cases leading to the
associahedron: points in convex position and the polytope of constrained
expansions in dimension $1$.
Section~\ref{towards} attempts to
   put the results in 1 and 2 dimensions into
   a broader perspective, with the aim of extending
   the results to higher dimensions and
   to point sets which are not in general position.
We conclude with some final comments in Section~\ref{final}.

\section{Preliminaries}
\label{preliminaries}

\paragraph{\bf Abbreviations and conventions.} Throughout this paper we will
assume general position for our point sets,
i.e. we assume that no $d+1$ points in $\reals^d$
lie in the same hyperplane (unless otherwise specified).
We abbreviate  ``polytope
of pointed pseudo-triangulations'' as \emph{ppt-polytope},
``one-degree-of-freedom mechanism'' as \emph{1DOF mechanism}
and ``pseudo-triangulation expansive mechanism'' as \emph{pte-mechanism}.

For an ordered sequence of $d+1$ points $q_0,\dots,q_d\in \reals^d$,
$\det(q_0,\dots,q_d)$ denotes 
the determinant of the $(d+1)\times(d+1)$ matrix with columns
$(q_0,1), \dots,\allowbreak (q_d,1)$. Equivalently, $\det(q_0,\dots,q_d)$ equals 
$d!$ times the Euclidean volume of the simplex 
with those $d+1$ vertices, with a sign
depending on the orientation.

\paragraph{\bf Pseudo-triangulations.}

A \emph{pseudo-triangle} is a simple polygon with only three convex vertices
(called \emph{corners})
joined by three inward convex polygonal chains,
see Figure~\figref{pseudot}a.
In particular, every triangle is a pseudo-triangle.
A \emph{pseudo-triangulation} is a partitioning
of the convex hull of a point set $P=\{p_1, \ldots , p_n\}$
into pseudo-triangles using $P$ as vertex set.

Pseudo-triangulations are \emph{graphs embedded on $P$},
i.e. graphs drawn in the plane on the vertex set $P$ and with
straight-line segments as edges.
We will work with other graphs embedded in the plane.
If edges intersect only at their end-points,
as is the case for pseudo-triangulations, the graphs
will be called \emph{non-crossing} or \emph{plane graphs}.
A graph is \emph{pointed at a vertex $v$} if there is (locally) an
angle at $v$ strictly larger than $\pi$ and containing no edges. Under our
general position assumption, convex-hull
vertices are pointed for any  
embedded graph, as are vertices of degree at most
two. A graph is called \emph{pointed} if it is pointed at every vertex.
Parts (b) and (c) of Figure~\figref{pseudot} (including the broken edges)
show two pointed pseudo-triangulations of a certain point set.

\begin{figure}[htb]
  \centering
{\def\IPEfile{pseudos.ipe}\input{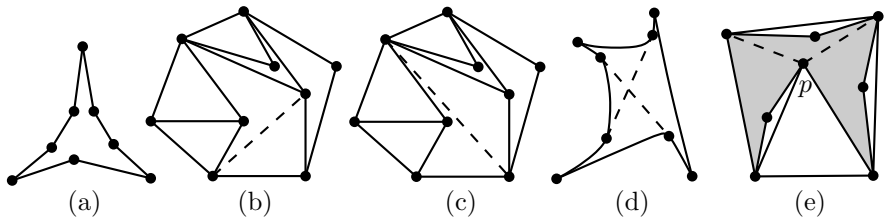}}
\caption{(a)
A pseudo-triangle.
(b) A minimum, or pointed, pseudo-triangulation.
(c) The broken edge in~(b) is flipped, and gives another
pointed pseudo-triangulation.
(d)~A schematic drawing of the flip operation.
(e)~The two edges involved in a flip may share a common vertex~$p$.
}
\figlab{pseudot}
\end{figure}

The following properties of pseudo-triangulations were
initially proved for
the slightly different situation of
pseudo-triangulations of \emph{convex objects} by
Pocchiola and Vegter~\cite{PV,PV2}.
For completeness, we sketch the easy proofs (see also~\cite{BKPS}).

\begin{lemma}
\label{lem:prop-pointed}
  {\bf (Streinu \cite{Streinu-2000})}
Let $P$ be a set of $n$ points in general position in the plane. Let $G$ be a
pointed and non-crossing graph on $P$.
\begin{itemize}
\item[\rm(a)] $G$ has at most $2n-3$ edges, with equality if and only if it is a
  pseudo-triangulation.
\item[\rm(b)] If $G$ is not a pseudo-triangulation, then edges can be added to it
  keeping it non-crossing and pointed. 
\end{itemize}
\end{lemma}

\begin{proof}
(a) 
If the graph $G$ is not connected, we can analyze the
components separately. So let us assume that it is connected.
Let $e$ and $f$ be the numbers of
edges and bounded faces in $G$. 
Let $a^+$ and $a^-$ denote the number of angles which
are $>\pi$ and $<\pi$, respectively. (By
our general position assumption, there are no angles equal to $\pi$.)
Clearly, $a^++a^-=2e$. Pointedness means $a^+=n$ and, since any bounded face
has at least three convex vertices, $a^-\ge 3f$ with equality if and only if
$G$ is a pseudo-triangulation.
The equation $2e \ge n+3f$, together with
Euler's formula $e=n+f-1$, implies $e\le 2n-3$ (and $f\le n-2$).

(b) The basic idea is that the
addition of geodesic
paths (i.e., paths which have shortest length among those sufficiently close
to them) between convex vertices of a polygon keeps the graph pointed and 
non-crossing and, unless the polygon is a pseudo-triangle, there is always some
of these geodesic paths going through its interior.
\end{proof}

 Streinu \cite{Streinu-2000} proved 
 the following additional properties of pointed
pseudo-triangulations, which we do not need for our results but
which may interest the reader:
\begin{itemize}
\item Every pseudo-triangulation on $n$ points has at least $2n-3$ edges, with
  equality if and only if it is pointed. Hence, pointed pseudo-triangulations
  are the pseudo-triangulations with the minimum number of edges. For this
  reason they are called {\em minimum} pseudo-triangulations in
  \cite{Streinu-2000}. In contrast with part (b) of Lemma
  \ref{lem:prop-pointed}, not every
  pseudo-triangulation contains a pointed one. An example of this is a regular
  pentagon with its central point, triangulated as a wheel. Hence, a {\em
  minimal} pseudo-triangulation is not always pointed.

\item The graph of any pointed pseudo-triangulation has the Laman property: it
  has $2n-3$ edges and the subgraph induced on any $k$ vertices has at most
 $2k-3$ edges. This property characterizes generically minimally rigid graphs
 in the plane (\cite{laman}, see also \cite{gss});
that is, graphs which are minimally rigid in almost all their
embeddings in the plane.

\item All pointed pseudo-triangulations can be obtained starting with a
  triangle and adding vertices one by one
and adding or adjusting edges, in
much the same way as the Henneberg construction of generically minimally rigid
graphs (cf.~\cite[page 113]{gss}), suitably modified to give pointed
pseudo-triangulations in intermediate steps (see the details in
\cite{Streinu-2000}). 
\end{itemize}

The other crucial properties of pointed
pseudo-triangulations that we use are that all interior edges
can be flipped in a natural way (part (a) of the following statement) and that
the graph of flips between pointed pseudo-triangulations of any point set
is connected. Both results were known to Pocchiola and Vegter for
{pseudo-triangulations of convex objects} (see \cite{PV,PV2}).
Parts (b) and (c) of Figure~\figref{pseudot} show a flip between pointed
pseudo-triangulations.
An $O(n^2)$ bound on the diameter of the flip graph is
proved in~\cite{BKPS}.

\begin{lemma}
\label{lem:props}
\label{lem:flips}
{\bf (Flips between pointed pseudo-triangulations)}
Let $P$ be a point set in general position in the plane.
\begin{enumerate}
\item[\textup{(a)}] \textup({\bf Definition of Flips}\textup)
\label{lem:unique}
When an interior edge \textup(not on the convex hull\textup)
is removed from a pointed
pseudo-tri\-angula\-tion of $P$,
there is a unique way to put back another edge
to obtain a different pointed pseudo-triangulation.

\item[\textup{(b)}]
 \textup({\bf Connectivity of the flip graph}\textup) 
 \label{lem:connect}
The graph whose vertices are pointed pseudo-triangulations and whose edges
correspond to flips of interior edges is connected.
\end{enumerate}
\end{lemma}

\begin{proof}
\cite{BKPS,Streinu-2000}
(a) When we remove an interior edge from a pointed
pseudo-triangulation we get a planar and pointed graph with $2n-4$ edges.
The
same arguments of the proof of Lemma~\ref{lem:prop-pointed} imply now that
$a^- = 3f + 1$. Hence, the new face created by the removal must be a
pseudo-quadrilateral (that is, a simple polygon with exactly four convex
vertices).

In any pseudo-quadrilateral there are exactly two ways of inserting
an interior edge to divide it into pseudo-triangles, which can be obtained by
the shortest paths between opposite convex vertices of the
pseudo-quadrilateral (see the details in Lemma 2.1 of \cite{Streinu-full}, and
a schematic drawing in
Figure~\figref{pseudot}d). One of these two is the
edge we have removed, so only the other one remains.
Note that the two interior edges of a 
 pseudo-quadrilateral may be incident to the same vertex,
see Figure~\figref{pseudot}e.
This can only happen when the interior angle at this vertex is
bigger than~$\pi$.

(b) Let $p$ be a convex hull vertex in $P$. 
Pointed pseudo-triangulations in which
$p$ is not incident to any interior edge are just pointed
  pseudo-triangulations of $P\setminus\{p\}$ together with
  the two tangent edges from $p$ to
  the convex hull of the rest. By
  induction, we assume all those pointed pseudo-triangulations
  to be connected to each other. To show that all others are also connected to
  those, just observe that
  if a pointed pseudo-triangulation has an interior edge
  incident to $p$, then a flip on that edge inserts an edge not incident to
  $p$.
(The case of Figure~\figref{pseudot}e cannot happen for a hull vertex~$p$.)
 Hence the number of interior edges incident to $p$ decreases.
\end{proof}

\paragraph {\bf Infinitesimal rigidity.}
In this paper we work mostly with points in dimensions $d=2$ and
$d=1$.
Occasionally we will use superscripts to denote the components of
the vectors 
$p_i=(p_i^1, \ldots, p_i^d)$.

An \emph{infinitesimal} \emph{motion} on a point set $P=\{p_1, \ldots,
p_n\}\in \reals^d$
is an assignment of a velocity vector
$v_i=(v_i^1, \ldots, v_i^d)$ to each  point $p_i$, $i=1, \ldots, n$.
The \emph{trivial infinitesimal motions} are those which
come from (infinitesimal) rigid transformations of the whole
ambient space. In $\reals^2$
these are the translations (for which all the $v_i$'s are equal vectors)
and rotations with a certain center $p_0$ (for which each $v_i$ is
perpendicular and proportional to the segment $p_0p_i$). Trivial motions form a
linear subspace of dimension $\binom{d+1}2$ in the linear space
$(\reals^d)^n$ of all infinitesimal motions.
 Two infinitesimal
motions whose difference is a trivial motion
will be considered equivalent, leading to a reduced
space of \emph{non-trivial} infinitesimal motions of
dimension $dn-\binom{d+1}{2}$. In particular,
this is $n-1$ for $d= 1$ and $2n-3$ for $d=2$.
Rather than performing a formal quotient of vector spaces
we will ``tie the framework down'' by fixing $\binom{d+1}{2}$ 
variables. E.g., for $d=1$ we can choose:
\begin{equation} \label{eq:normal1}
 v_1 = 0
\end{equation}
and for $d=2$ (assuming w.l.o.g. that $p_2^2 \ne p_1^2$):
\begin{equation} \label{eq:normal2}
 v_1^1 = v_1^2 = v_2^1 = 0
\end{equation}
Here, $p_1$ and $p_2$ can be any two vertices. A different choice
of normalizing conditions only 
amounts to a linear transformation in the space of
infinitesimal motions.

In rigidity theory, a graph $G=(P,E)$ embedded on $P$
is customarily called a \emph{framework} and denoted
by $G(P)$.
We will use the term framework when we want to emphasize
its rigidity-theoretic properties (stresses, motions), but we
will use the term graph when speaking about graph-theoretic
properties, even if graph is embedded on a set~$P$.
For a given framework $G=(P,E)$,
an infinitesimal motion such that $\langle p_i-p_j, v_i-v_j\rangle = 0$ for
every edge $ij\in E$ is called a \emph{flex} of $G$.
This condition states that the length of the edge $ij$ remains
unchanged, to first order.
The trivial motions are the flexes of the complete graph,
provided that the vertices span the whole space $\mathbb{R}^d$.
A framework is \emph{infinitesimally rigid} if it has
no non-trivial flexes.  It is \emph{infinitesimally flexible} or an
\emph{infinitesimal mechanism} otherwise.

Infinitesimal motions are
to be distinguished from global motions, which describe
paths for each point throughout some time interval.
In this paper we are not concerned with global motions, nor their associated
concept of rigidity, weaker than infinitesimal 
rigidity (\cite[Theorem~4.3.1]{Connelly-Whiteley-1996} or \cite[page 6]{gss}).
Let us also 
insist that we distinguish between {\em infinitesimal motions} (of the
point set) and {\em flexes} (of the framework or embedded graph), 
while the terms flex and
infinitesimal motion are sometimes equivalent in the rigidity theory
literature.

The (\emph{infinitesimal})
\emph{rigidity map} $M_{G(P)}\colon (\reals^d)^n \to \reals^{E(G)}$
is a linear map associated with an embedded framework $G(P)$,
$P\subset \reals^d$. It sends each infinitesimal motion
$(v_1,\dots,v_n)\in (\reals^d)^n$ to the vector of infinitesimal edge increases
$(\langle p_i-p_j, v_i-v_j\rangle)_{ij\in E}$.
When no confusion arises, it will be simply denoted as $M$.
The number $\langle p_i-p_j, v_i-v_j\rangle$ is called
the \emph{strain} on the edge $ij$ in the engineering literature.
As usual, the image of $M$ is denoted by
$\im M=\{\, f \mid f=Mv \,\}$.
The matrix of $M$ is called the rigidity matrix.
 % (also denoted by $M$, if no confusion arises),
In this matrix,
the row indexed by the edge $ij\in E$ has $0$ entries everywhere except in the
$i$-th and $j$-th group of $d$ columns, where the entries are $p_i-p_j$
and $p_j-p_i$, respectively.

The kernel of $M$ (after reducing $\reals^{dn}$ to $\reals^{dn-
\binom{d+1}{2}}$
  by forgetting trivial motions) is the space of flexes of $G(P)$.
In particular, a framework is infinitesimally 
rigid if and only if the kernel of its associated rigidity map $M$
is the subspace of trivial motions. In general, the dimension of the 
(reduced) space of
flexes is the  \emph{degree of freedom} (DOF) of the framework.
A 1DOF mechanism is a mechanism with {one degree of freedom}.

Finally, \emph{expansive (infinitesimal) motions}
$v_1,\ldots, v_n$ are those which
simultaneously increase (perhaps not strictly) all distances:
$\langle p_i-p_j, v_i-v_j\rangle \geq 0$ for every pair $i,j$ of
vertices. A mechanism is
\emph{expansive} if it has non-trivial expansive flexes.

The following results of Streinu~\cite{Streinu-2000} can be obtained as a corollary of our main
result (see the proof after the statement of Theorem \ref{main}).

\begin{proposition}\label{lem:pseudot}
{\bf (Rigidity 
 of pointed pseudo-trian\-gu\-lations~\cite{Streinu-2000})}
\begin{enumerate}
\item[\rm(a)]
Pointed pseudo-triangulations are  minimally infinitesimally rigid \textup(and
therefore rigid\textup).

\item[\rm(b)]
 The removal of a convex hull edge from a pointed pseudo-triangulation
yields a 1DOF expansive mechanism \textup(called a
\emph{pseudo-triangulation expansive mechanism}
or shortly a \emph{pte-mechanism}\textup).
\qed
\end{enumerate}
\end{proposition}
Part~(a) is in accordance with the fact that
the graph of any pointed pseudo-triangulation
has the Laman property, and hence is generically
rigid in the plane. It is a trivial consequence of~(a)
that the removal of an edge creates a (not necessarily expansive)
1DOF mechanism.
The expansiveness of pte-mechanisms (part~(b))
was proved in~\cite{Streinu-2000}
using the Maxwell-Cremona correspondence between self-stresses
and 3-d liftings of planar frameworks,
a technique that was introduced in~\cite{CDR}.

\paragraph {\bf Self-stresses.}
A \emph{self-stress} (or an \emph{equilibrium stress}) on a framework $G(P)$
(see \cite{Whiteley} or \cite[Section~3.1]{CDR})
is an assignment of scalars $w_{ij}$
to edges such that $\forall i\in P$, $\sum_{ij\in E} w_{ij}(p_i-p_j)=0$.
That is, the self-stresses are the row dependences of
the rigidity matrix~$M$.
The proof of the following lemma is then straightforward.

%In the following Lemma, whose proof is straightforward,
%$M^T\colon \reals^{E(G)}\to (\reals^d)^n$ denotes the
%adjoint
%linear map of $M$, whose matrix is the transpose of the matrix of $M$.

\begin{lemma}
\label{lem:stress4}
Self-stresses form the orthogonal complement of the linear subspace 
$\im M\subset \reals^{\binom d 2}$.  In other words,
$(w_{ij})_{ij\in E}$ is a self-stress if and only if for every infinitesimal
motion $(v_1,\dots,v_n)\in (\reals^d)^n$ the following identity holds:
$$
\sum_{ij\in E}w_{ij}\langle p_i-p_j, v_i-v_j\rangle = 0
%\eqno 
\qed
$$
\end{lemma}

As an example, the following result gives explicitly a stress
for the complete graph on any affinely dependent point set:

\begin{lemma}
\label{lem:dependent}
Let $\sum_{i=1}^n \alpha_i p_i=0$, $\sum \alpha_i=0$,  be an affine dependence
on a point set $P=\{p_1,\dots,p_n\}$. Then, $w_{ij}=\alpha_i\alpha_j$
for every $i,j$ defines a self-stress of the complete graph $G$ on $P$.
\end{lemma}

\proof % \begin{proof}
For any  $p_i\in P$ we have:
\[
\sum_{ij\in G} w_{ij}(p_i-p_j)=
\sum_{j=1}^n \alpha_i \alpha_j(p_i-p_j)=
\alpha_i p_i  \sum_{j=1}^n  \alpha_j
- \alpha_i \sum_{j=1}^n \alpha_j p_j = 0.
\qed
\]
% \end{proof}

Let us analyze here the case of $d+2$ points $P=\{p_1,\dots,p_{d+2}\}$
in general position in $\reals^d$ (this is the first non-trivial case,
because no self-stress can arise
between affinely independent points).
It can be easily checked that, under these assumptions,
removing any single edge from the
complete graph on $P$ leaves a minimally infinitesimally rigid graph. This
implies that the complete graph has a unique self-stress (up to a scalar
factor). This self-stress is the one given in Lemma~\ref{lem:dependent} for
the unique % (up to  a scalar factor) 
affine dependence on $P$. The coefficients
of this dependence can be written as:
\[
\alpha_i = (-1)^{i}\det([p_1,\dots,p_{d+2}]\backslash \{p_i\}).
\]
(Recall that $\det(q_0,\dots,q_d)$ is
$d!$ times the signed volume of the simplex spanned by the
$d+1$ points $q_0,\dots,q_d\in \reals^d$.)

The special case $d=2$, $n=4$ will be extremely 
relevant to our purposes, and it will be convenient 
to renormalize the unique self-stress as follows:

\begin{lemma}
\label{lem:stresses}
The following gives a self-stresses for any four points $p_1,\dots,p_4$ in
general position in the plane\textup:
\begin{equation}
\label{eq:stress}
w_{ij}:=\frac{1} {\det(p_i,p_j,p_k)  \det(p_i,p_j,p_l) }
\end{equation}
where $k$ and $l$ are the two indices other than $i$ and $j$.
\end{lemma}

\proof % \begin{proof}
Set $\alpha_i=(-1)^{i}\det([p_1,\dots,p_4]\backslash\{p_i\})$
in Lemma~\ref{lem:dependent}
and divide all the $w_{ij}$'s
of the resulting self-stress by the non-zero constant
\[
-\det(p_1,p_2,p_3)\det(p_1,p_2,p_4)\det(p_1,p_3,p_4)\det(p_2,p_3,p_4).
%\eqno 
\qed
\]
%\end{proof}

The direct application of Lemma~\ref{lem:dependent} would give as
$w_{ij}$ a product of two determinants, rather than the inverse of such a
product. The reason why we prefer 
the self-stress of Lemma~\ref{lem:stresses} 
is because of the signs it produces. The reader
can easily check, considering the two cases of four points in convex
position and one point inside the triangle formed by the other three, that
with the choice of Lemma~\ref{lem:stresses} boundary edges always receive
positive stress and interior edges negative stress. This uniformity is good
for us because in both cases pointed pseudo-triangulations are the graphs
obtained deleting from the complete graph any single interior edge.

%The self-stress of four points is unique up
%to a constant factor. Multiplying the formulas
%(\ref{eq:stress}) by the product of the four determinants of the four
%points we would get somewhat simpler formulas, each $w_{ij}$ becoming the
%product of two determinants, 
%Our choice is more convenient because no matter
%which four points we take, the $w_{ij}$'s will be positive for
%convex hull edges and negative for interior edges.
%Note that a framework with five edges on four points is non-crossing and
%pointed (i.e. a pointed pseudo-triangulation) if and only if the
%missing edge is interior. We will use these observations in our construction
%of the ppt-polytope.

\paragraph{\bf The expansion cone.}
We are given a set of $n$ points
$P=(p_1, \ldots, p_n)$ in $\reals^d$
that are to move with (unknown) velocities $v_i\in\reals^d$, $i=1,\ldots,n$.
An \emph{expansive} motion is a motion in which no inter-point distance
decreases. This is described by the system of homogeneous linear
inequalities:
\begin{equation} \label{eq:0}
\langle p_i-p_j, v_j-v_i \rangle \ge 0,\qquad \forall\ 1\le i<j\le n
\end{equation}
and hence defines a polyhedral cone.

The only motions in the intersection of all facets of the cone are the
trivial motions.
Thus, when we add normalizing equations like \eeqref{normal1} or~\eeqref{normal2},
we get a \emph{pointed polyhedral cone}
containing the origin as a vertex.
% , in the reduced space of infinitesimal motions
We call it the cone of expansive motions
or simply the \emph{expansion cone} of $P$.

An extreme ray of the
expansion cone is given by a  maximal set of inequalities
satisfied with equality by non-trivial motions.
Each inequality corresponds to an edge of the point set, so that the ray
corresponds to a graph embedded in our point set.
The cardinality of this set of edges is at least
the dimension of the cone minus 1,
but may be much larger. Let's analyze the low dimensional cases.

For $d=1$ the expansion cone is not very interesting. Let's assume that the
points $p_i\in \reals$ are labeled in increasing order $p_1<p_2<\cdots <
p_n$. Then:

\begin{proposition} \lemlab{X1}
The expansion cone in one dimension has $n-1$ extreme rays corresponding to
the motions where $p_1,\ldots,p_i$ remain stationary and the points
$p_{i+1},\ldots,p_n$ move away from them at uniform speed\textup:
\begin{equation}
\label{eq:exp1}
0=v_1=v_2=\cdots=v_i < v_{i+1}= \cdots = v_n
\end{equation}
\end{proposition}

\begin{proof}
Note that the actual values of $p_i$ are immaterial in this case.
The expansion cone is given by the linear system $v_j\geq v_i$, $1\leq i\leq j
\leq n$ plus the extra condition $v_1=0$, and any maximal set of inequalities
satisfied with equality and yet not trivial is obviously given by
(\ref{eq:exp1}).
\end{proof}

The $2$d case is more complex and requires additional terminology.
% Generic 
1DOF mechanisms may contain rigid subcomponents (r-components,
cf.~\cite{gss}): maximal sets of some $k$ vertices spanning a Laman subgraph
on $2k-3$ edges.  The r-components of pte-mechanisms are
themselves pseudo-triangulations spanning convex subpolygons
including all points in their interior.
%We say that
%two pt-mechanisms are \emph{r-equivalent} if they differ only inside a rigid
%subcomponent.
Adding edges to complete each r-component to a complete subgraph yields a
\emph{collapsed pte-mechanism} (see Figures~\figref{collapse}
and~\figref{extreme}).

%The equivalence classes will be called \emph{pt-mechanisms with
%collapsed rigid subcomponents} or simply (when no confusion arises as to the
%removed edge) a \emph{collapsed pte-mechanism}.

\begin{figure}[htb]
  \centering
\includegraphics% [width=0.8 \textwidth, height=0.423\textwidth]
{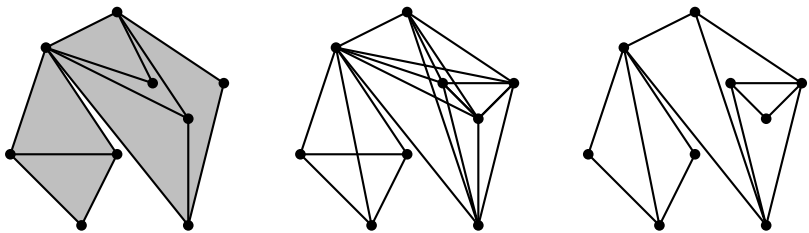}
\caption{A pte-mechanism
with rigid sub-components (convex subpolygons) drawn shaded,
the corresponding
collapsed pte-mechanism, and another
pte-mechanism that yields the same expansive motion.
}
\figlab{collapse}
\end{figure}

\begin{figure}[htb]
  \centering
\includegraphics[width=0.8 \textwidth, height=0.423\textwidth]
{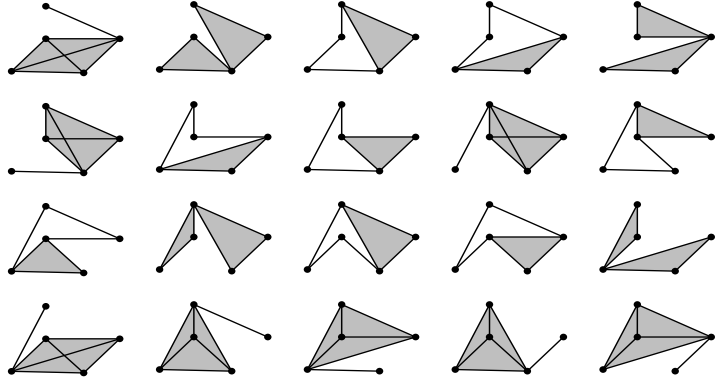}
\caption{
The collapsed pte-mechanisms corresponding to the 20 extreme rays of the
expansion cone for a planar point set of $5$~points.
The rigid sub-components (complete subgraphs) are shaded.
}
\figlab{extreme}
\end{figure}

%\begin{figure}[htb]
%  \centering
%%  \includegraphics[width=1.0 \textwidth, height=0.5\textwidth]{extreme.ai}
%%\psfig{figure=extreme.ai,height=1.7in,width=2.8in}
%\caption{
%The collapsed pte-mechanisms corresponding to the extreme rays of the
%expansion cone for a planar point set on $5$ points.
%}
%\figlab{extreme}
%\end{figure}

\begin{proposition}
\lemlab{ex-rays}
In dimension $2$, the extreme rays of the expansion cone 
correspond to the collapsed pte-mechanisms.
%equivalence classes of pt-mechanisms modulo r-equivalence.
\qed
\end{proposition}

The proof will be given in Section~\secref{x-cone},
after we have determined the extreme rays of a perturbed version of
the polytope.
%We will not use this statement in the rest of the paper.
%\begin{proof}
%We sketch a proof along the lines
% of the proof
%technique used in \cite{CDR} and \cite{Streinu-2000},
%using the Maxwell-Cremona correspondence between self-stresses
%and 3-d lifting of planar frameworks.
%Another independent proof will be given in Section~\secref{x-cone}.

%One convex hull
%edge must be missing, because otherwise the corresponding motion won't be
%expansive (cf.~Lemma~\lemref{bounded} below).
%An extreme ray is one-dimensional, hence the corresponding set of
%edges induces a 1DOF mechanism. The r-components of the mechanism must have all
%the edges present in the extreme ray, since they are linearly dependent.
%Thus to show that pte-mechanisms correspond to extreme rays we must show that
%the mechanism is expansive, which is done via a direct adaptation of the proof
%technique from \cite{CDR} and \cite{Streinu-2000} (by contradiction, using
%self-stresses,
%$3$-d liftings and Maxwell's theorem).

%To complete the proof we must show that the graph corresponding to an extreme
%ray has no crossing edges and no not-pointed vertices outside the (interior of
%the) r-components. Suppose it has two crossing edges that move relative to each
%other in the mechanism. Then their extreme vertices can't simultaneously move
%away from each other, contradicting the expansive property and the assumption
%that we started with an extreme ray of the expansion cone. Similarly, if a
%vertex is not pointed but its adjacent edges move relative to each other, then
%the overall motion can't be expansive.
%\end{proof}

\paragraph{\bf The polytope of constrained expansions.}
The construction we will give in section \ref{polytope}
can roughly be interpreted as separating the
pseudo-triangulations contained in the same collapsed
pte-mechanisms, to obtain a polyhedron whose vertices correspond to distinct
pseudo-triangulations.  The original expansion cone is highly degenerate: its
extreme rays contain information about \emph{all} the bars whose length is
unchanged by a motion of a 1DOF expansive mechanism.  We would like to
perturb the constraints~(\ref{eq:0}) to eliminate these degeneracies and
recover pure pseudo-triangulations. We do so by giving up homogeneity, i.e.,
by translating the facets of the expansion cone. Our system will become:
\begin{equation} \label{eq:f}
  \langle v_j-v_i, p_i-p_j \rangle  \ge f_{ij},
  \qquad\ \forall 1\le i<j\le n
\end{equation}
for some numbers $f_{ij}$.
In some cases we will change these inequalities to equations
for the edges on the convex hull of the given point set.
\begin{equation} \label{eq:f=}
  \langle v_j-v_i, p_i-p_j \rangle  = f_{ij},
  \qquad\ \text{for the convex hull edges $ij$}.
\end{equation}

Section \ref{polytope} proves our main result, Theorem~\ref{main}:
For any point set in general position in the plane and for
some appropriate choices of the parameters $f_{ij}$, \eeqref{f}
defines a polyhedron whose
vertices are in bijection with pointed pseudo-triangulations and all
lie in a unique maximal bounded face given by~\eeqref{f=}.
We call this face
the ``polytope of pointed
pseudo-triangulations'' or ppt-polytope. 

A similar thing in 1d is done in
Section \ref{1dim}, with the surprising outcome that the (unique) maximal
bounded face of the polyhedron turns out to be an associahedron with
vertices corresponding to \emph{non-crossing alternating trees}
(which are Catalan structures, as shown in \cite{GGP}). The next paragraph
prepares the ground for this result.

\paragraph{\bf The associahedron.}
The associahedron is a polytope which has a vertex for
every triangulation of a convex $n$-gon, and in which two vertices are
connected by an edge of the polytope
if the two triangulations are connected by an edge flip. Equivalently, various
types of Catalan structures are reflected in the associahedron.
Fig.~\ref{fig:assoc} shows an example.

\begin{figure}[htb]
  \centering
%\includegraphics[width=0.35 \textwidth, height=0.40\textwidth]{bassoc.ipe}
%%\psfig{figure=bassoc.ipe,height=1.7in,width=3.0in}
%\ \ \ \ \ \ \ \ \ \ \ \ \ \ \ \ \ \ \ 
{\def\IPEfile{assocs.ipe}\input{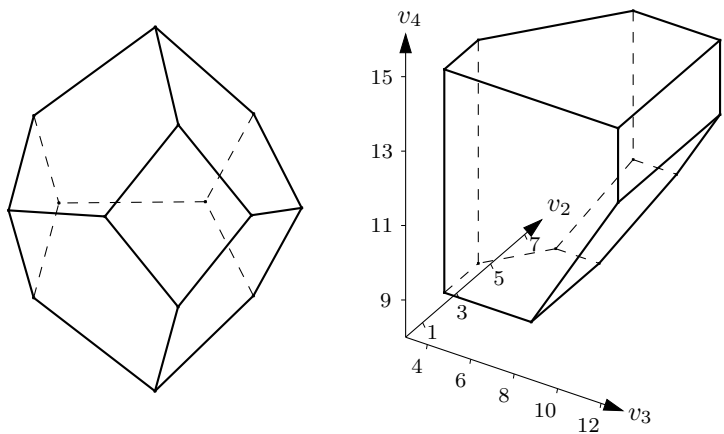}}
    \caption{The three-dimensional associahedron.
The vertices represent all triangulations of a convex hexagon or
all possible ways to insert parentheses into the product $a*b*c*d*e$.
Left: a symmetric representation,
the secondary polytope of a regular hexagon.
Right: our representation, from Section~\ref{1dim}.
Both pictures are orthogonal projections.
}
\label{fig:assoc}
\end{figure}

There is an easy geometric realization of this polytope associated with each set
$P$ of $n$ points \emph{in convex position} in the plane,
as a special case of a secondary polytope
(Gel$'$fand, Zelevinski\v\i, and Kapranov~\cite{GZK}, see
also Ziegler~\cite[Section~9.2]{Ziegler}).
Every triangulation is represented by a vector $(a_1,\ldots,a_n)$ of
$n$ component with the entry $a_i$ being simply the sum of the areas
of all triangles of
the triangulation that are incident to the $i$-th vertex.
We will refer to this realization as the
\emph{classical realization} of the associahedron.
It %This realization of the associahedron
depends on the location
of the vertices of the convex $n$-gon, but all polytopes
that one gets in this way are combinatorially equivalent.
Their face lattice is the poset of polygonal subdivisions of the $n$-gon or,
in the terminology of the previous paragraphs, non-crossing and pointed graphs
embedded in $P$ and containing the $n$ convex hull edges. But observe that
the word ``pointed'' is superfluous for a graph with vertices in convex
position. The order
structure in this poset is just inclusion of edge sets (in reverse order since
maximal graphs represent vertices).

Dantzig, Hoffman, and Hu~\cite[Section~2]{DHH}, and independently
de Loera et al.~%
\cite{LHSS} in a more general setting,
have given other representations of the set of triangulations
as the vertices of a 0-1-polytope in $\binom n3$ variables
corresponding to the possible triangles of a triangulation
(the \emph{universal} polytope), or in $\binom n2$ variables
corresponding to the possible edges of a triangulation.
These realizations are in a sense most natural, but they have
higher dimensions and have more adjacencies between vertices
than the associahedron. Every classical associahedron, however, arises
as a projection of the universal polytope.
The first published realization of an associahedron is due
to Lee~\cite{Lee}, but it is not fully explicit.
A few earlier and more complicated ad-hoc
realizations that were never published
are mentioned in Ziegler~\cite[Section~0.10]{Ziegler}.

In this paper the associahedron appears in two forms. First, we will show that
for $n$ points in convex position our polytope of pointed
pseudo-triangulations is affinely equivalent to the secondary polytope of the
configuration, which is a classical associahedron (Section
\ref{convex2dim}). Second, as mentioned before, our construction adapted to
a one-dimensional point configuration produces in a natural way
an associahedron (Section \ref{1dim}).
Notice that in dimension $1$ the coordinates $p_i$ can be eliminated from
the constraints~(\ref{eq:f}).
% complaint{Changed f to g. GR}
Only the order of points along the line
matters. One can also look at the whole \emph{arrangement} of hyperplanes
of the form
\begin{equation} \nonumber % \label{eq:f}
  v_j-v_i = g_{ij}.
\end{equation}
This arrangement, for various special values of $g$,
has been the object of extensive combinatorial studies.
For $g\equiv0$ it is the classical \emph{Coxeter} or \emph{reflection}
arrangement of type $A_n$. The case $g\equiv1$ has been studied
by Postnikov and Stanley~\cite{PS}. The expansion cone of a 1d point set
is the \emph{positive cell} in the
arrangement $A_n$, and our associahedron is a bounded face of the polyhedron
obtained by translating the facets of this cell. 

A different realization of the associahedron
based in the root system of type $A_n$ has recently appeared in \cite{ChFoZe}. It is
interesting that these two new associahedra are not affinely equivalent to any
classical associahedron obtained as a secondary polytope, or to one another. 
Also, that we are trying to get a simple polyhedron, in
contrast to the above-mentioned choices of $g$ which lead
to highly degenerate arrangements.

\section{The Main Result: the Polytope of Pointed Pseudo-Triangulations}
\label{polytope}
In this section we prove our main result.
\begin{theorem}
\label{main}
For every set $P=\{p_1,\dots,p_n\}$
of $n\ge3$ planar points in general position, 
there is a choice of $f_{ij}$'s for which equations~\eeqref{f} 
together with the normalizing equations~\eeqref{normal2}
define a simple polyhedron $\phce_f(P)$ of dimension $2n-3$
with the following properties\textup:
\begin{enumerate}
\item The face poset of the polyhedron equals the opposite of the
  poset of pointed and non-crossing graphs on $P$, by the map sending each
  face to the set of edges 
  whose corresponding equations~\eeqref{f} are satisfied with equality over
  that face.
  In particular:
\begin{enumerate}
  \item[\textup{(a)}]
 Vertices of the polyhedron are in 1-to-1 correspondence with 
  pointed pseudo-triangulations of~$P$.
  \item[\textup{(b)}] Bounded edges correspond to flips of interior
  edges in pseudo-triangulations, i.e., to pseudo-triangulations with one
  interior edge removed.
  \item[\textup{(c)}] \label{ex-rays}
    Extreme rays correspond to pseudo-triangulations with one convex hull
    edge removed.
\end{enumerate}
\item The face $\pce_f(P)$ obtained by changing to
 equalities \eeqref{f=} those inequalities from \eeqref{f} which correspond 
 to convex hull edges of $P$ is bounded \textup(hence a polytope\textup)
 and contains all vertices. In other words, it is the unique
 maximal bounded face, and its 1-skeleton is the graph of flips among pointed
 pseudo-triangulations.
\end{enumerate}
\end{theorem}

The proof is a consequence of lemmas proved throughout this section.
Theorem \ref{thm:valid-new} states that the choice
$f_{ij} := \det(a,p_i,p_j)\det(b,p_i,p_j)$
produces the desired object, where $a$ and $b$ are 
any fixed points in the plane.
In Section \ref{meta-ppt} we will 
derive a more canonical description of the
polyhedron (and the polytope) in question.

It turns out that the polyhedron $\phce_f(P)$ is the most convenient object
for the proof.
The properties of the polytope $\pce_f(P)$ (part~2 of the theorem)
and the extreme ray description of the expansion cone
(Proposition~\lemref{ex-rays}), which may be more interesting
by themselves, are then easily derived.

Before going on, 
let us see that Theorem~\ref{main} implies
Proposition~\ref{lem:pseudot}. Observe that a framework is minimally
infinitesimally rigid if and only if the hyperplanes 
$\langle p_i-p_j, v_i-v_j\rangle=0$ corresponding to its edges $ij$
meet transversally and at a single point, in the
$(2n-3)$-dimensional space given by equations~\eeqref{normal2}. 
Part 1.a of our theorem says that this happens for the $2n-3$ translated
hyperplanes $\langle p_i-p_j, v_i-v_j\rangle=f_{ij}$ corresponding
to any pointed pseudo-triangulation, hence giving part~(a) of 
Proposition~\ref{lem:pseudot}. An (infinitesimally) expansive 1DOF mechanism 
is one whose corresponding hyperplanes intersect in a line contained in the
expansion cone. Part 1.c of the theorem says that this happens for a
pointed pseudo-triangulation with one hull edge removed, giving 
part~(b) of
Proposition~\ref{lem:pseudot}.

%The solutions $v\in {(\reals^{2})}^n$ of the perturbed system (\ref{eq:f})
%$Mv \geq f$
%define the 
%\emph{polyhedron of constrained expansions},
%which will
%be shown to give the desired object
%(the ppt-polyhedron) for \emph{certain} choices of $f$'s.
%If we change the inequalities for the convex hull edges to equations,
%we get a bounded face of the polyhedron,
%which will then be the 
%\emph{polytope of constrained expansions},
%(the ppt-polytope).

\paragraph{\bf The polyhedron and the polytope of constrained expansions.}
% {\bf Tight edges.}
The solution set $v\in {(\reals^{2})}^n$ 
of the system of inequalities~\eeqref{f}
together with the normalizing equations~\eeqref{normal2}
will be called the \emph{polyhedron of constrained expansions} $\phce_{f}(P)$ 
for the set of points $P$ and perturbation parameters ({constraints})
$f$.
We will frequently omit the point set $P$ when it is clear from the context.
A solution $v$
may satisfy some of the inequalities in (\ref{eq:f}) with equality: the
corresponding edges $E(v)$ of $G$ are said to be \emph{tight} for that
solution. In the same way, for a face $K$ of $\phce_{f}$ we call 
{\it tight edges}
of $K$ and denote $E(K)$ the edges whose equations are satisfied with equality
over $K$ (equivalently, over a relative interior point of $K$).
This is the correspondence that
Theorem~\ref{main} refers to: the edges $E(K)$ of a face $K$ form the
pointed and non-crossing graph corresponding to that face.

When $f\equiv 0$, we just get the expansion cone $\xcone$ itself
(in this sense, our notations $\phce_{0}$ and $\phce_{f}$ are consistent.)
This cone equals the recession cone of $\phce_f$, for any choice of~$f$.
(The recession cone of a polyhedron is the cone of vectors
parallel to infinite rays contained in the polyhedron.)
We will first establish a few properties of the expansion cone.

\begin{lemma} \lemlab{xcone-tight}
  \begin{enumerate}
  \item[\rm{(a)}]
The expansion cone $\phce_{0}$ is a pointed polyhedral
cone of full dimension $2n-3$ in
the subspace defined by the three
equations~\eeqref{normal2}.
\textup(In this context,
``pointed'' means that the origin is a vertex of the cone.\textup)
  \item[\rm(b)]
Consider the set $E(v)$ of tight edges for any feasible point
$v\in\xcone$. If $E(v)$ contains
\begin{enumerate}
\item [\textup{(i)}]
two crossing edges,
\item [\textup{(ii)}]
a set of edges incident to a common vertex with no angle larger
than~$\pi$ \textup(witnessing that $E(v)$ is not
pointed at this vertex\textup),
or
\item [\textup{(iii)}]
a convex subpolygon,
\end{enumerate}
then $E(v)$ must contain the complete graph between the endpoints
of all involved edges.
In case~\textup{(iii)}, this complete graph also includes all points
inside the convex subpolygon.
  \end{enumerate}
\end{lemma}
\begin{proof}
(a)~The dilation (scaling motion)
$v_i := p_i$ satisfies all inequalities~\eeqref{0} strictly.
By adding a suitable rigid motion,
the three equations~\eeqref{normal2} can be satisfied, too,
without changing the status of the inequalities~\eeqref{0},
and so we get a relative interior point in the
$(2n-3)$-dimensional subspace~\eeqref{normal2}.

If the cone were not pointed, it would contain two opposite
vectors $v$ and $-v$. From this we would conclude that
$\langle v_j-v_i,p_j-p_i\rangle = 0$ for all $i,j$, and hence
$v$ would be a flex of the complete graph on $P$. By the
 normalizing equations~\eeqref{normal2}, $v$ must then be~$0$.

(b)~We first consider~(iii),
which is the most involved case.
Let $v$ be an expansive motion which preserves all edge lengths of
some convex polygon.
First we see
that $v$ preserves all distances between polygon vertices:
indeed, if it preserves lengths of polygon edges but is not a trivial
motion of the polygon then the angle at some
polygon vertex $p_i$ infinitesimally decreases, because the sum of
angles remains constant.
But decreasing the angle at $p_i$ while preserving the lengths of the 
two incident edges
implies that the distance between the two vertices adjacent to $p_i$
in the polygon decreases. This is a contradiction.

By choosing $p_1$ and $p_2$ in~\eeqref{normal2}
to be polygon vertices, the above implies that
the polygon remains stationary under~$v$.
Now no interior point $p_i$ can
move with respect to the polygon,
without decreasing the distance to some polygon vertex:
If $v_i\ne 0$, there is at least
one hull vertex $p_j$ in the half-plane $\langle p_i-p_j,v_i\rangle < 0$.
The edge $ij$ will then violate condition~\eeqref{0}.

Case~(ii) is similar:
If the edges incident to a vertex $p_i$ do not move rigidly,
at least one angle between two neighboring edges must decrease,
and, this angle being less than~$\pi$, this implies that
the distance between the endpoints of these edges decreases,
a contradiction.

For case~(i), we apply
Lemma~\ref{lem:stress4}
to our given four-point set
in convex position,
 with the self-stress of
Lemma~\ref{lem:stresses},
which is positive for the four hull edges and negative
for the two diagonals.
This implies that this four-point set
can have no non-trivial expansive motion which is
not strictly expansive on at least one of the two diagonals.
% will be proved later as a corollary of
% Lemma~\ref{lem:incomplete}.
\end{proof}

As an immediate consequence of Lemma~\lemref{xcone-tight}(a),
we get:
\begin{corollary} \lemlab{nonempty}
  $\phce_f(P)$ is a $(2n-3)$-dimensional unbounded polyhedron with at least
  one vertex, for any choice of parameters $f$.
\qed
\end{corollary}

It is easy to derive part 2 of Theorem~\ref{main} from part 1.
For every vertex or bounded edge of $\phce_f(P)$, the set
$E(v)$ contains all convex hull edges of $P$. On the contrary, for
any unbounded edge (ray) of $\phce_f(P)$, the set
$E(v)$ misses some convex hull edge of $P$.
Hence, by setting to equalities
the inequalities corresponding to convex hull edges
we get a face $\pce_{f}(P)$ of $\phce_{f}(P)$
which contains all vertices and bounded edges of $\phce_f(P)$, but 
no unbounded edge.

In order to prove part~1, we first need to check that indeed $\pce_{f}$
is a bounded face, and hence a polytope which we call
the \emph{polytope of constrained expansions}
or pce-polytope
for the set of points $P$ and perturbation parameters~$f$.

\begin{lemma}
  \lemlab{bounded}
 For any choice of~$f$, 
$\pce_f(P)$ is a bounded set.
\end{lemma}
\begin{proof}
Suppose that $v_0+tv$ is in $\pce_{f}$ for all $t\ge 0$. Then
we must have $v\in\pce_{0}$. Hence, it suffices to show that
$\pce_{0}=0$, i.e. that the framework consisting of all convex hull
edges has no non-trivial expansive flexes.
This is an immediate consequence of
 Lemma~\lemref{xcone-tight}b(iii).
\end{proof}

%Since we are only interested in graphs containing all the convex hull
%edges, we set all equalities for convex hull edges to equalities in
% our system 
%(\ref{eq:f}). This has the consequence that our polytope is now bounded, 
%and will be called the 
%\emph{polytope of constrained expansions} $\pce_{f}$ 
%for the set of points $P$ and perturbation parameters (\emph{constraints}) $f$. 
%\compaint{we can perhaps use ``constraints" instead of ``perturbation
%parameters". Paco.} 
%\compaint{I am not sure, let's see what G\"unter says, constraints may mean 
%the whole inequalities, not just the constants on the right hand side. Ileana}

\paragraph{\bf Reducing the problem to four points.}
We call a choice of the constants $f=(f_{ij})\in \reals^{\binom n 2}$ 
\emph{valid} if the corresponding polyhedron $\phce_{f}$
of constrained expansions 
has the combinatorial structure claimed in Theorem~\ref{main}.

\begin{lemma}
\label{lem:planar}
A choice of $f\in \reals^{\binom n 2}$ is valid 
if and only if the graph $E(v)$ of tight edges
corresponding to any feasible point
$v\in\phce_f(P)$ 
is non-crossing and pointed.
\end{lemma}

\begin{proof}
Necessity is trivial, by definition of being valid.
To see sufficiency note that,
by Corollary~\lemref{nonempty}, $\phce_{f}$ has dimension $2n-3$.
Thus, any vertex $v$ of the polyhedron is incident to
at least $2n-3$ faces $E(v)$.
If $E(v)$ is non-crossing and pointed,
%Lemma~\ref{lem:properties}(\lemref{characterization})
Lemma~\ref{lem:prop-pointed}
implies that $E(v)$ has \emph{exactly} $2n-3$ incident
faces and is a pointed pseudo-triangulation. 
In particular, the polyhedron is
simple. Also, since the tight edges of faces incident to $v$ are different
subgraphs of $E(v)$, the poset of faces incident to the vertex $v$
is the poset of all subgraphs of the pointed pseudo-triangulation $E(v)$.

It remains only to show
that every pointed pseudo-triangulation actually appears
as a vertex, for which we use a somewhat
indirect argument, based on the fact that the flip graph is connected.
This type of argument
has also been used by Carl Lee
for the case of the associahedron in~\cite{Lee},
where it is attributed to Gil Kalai and Micha Perles.

Since the polytope is simple, every vertex $v$ is incident to $2n-3$ edges
of $\phce_{f}$. The sets of tight edges corresponding to them are the $2n-3$
subgraphs of $E(v)$ obtained removing a single edge. We denote by $K_{ij}$ the
polyhedral edge corresponding to the removal of $ij$.
By Lemma~\lemref{bounded}, if $ij$ is interior then $K_{ij}$ is
bounded. Hence, it is incident to another vertex, which must
correspond to a pointed pseudo-triangulation that completes
$E(v)-\{ij\}$.
By Lemma~\ref{lem:props}(a),
this can only be the one obtained from $E(v)$ by a flip at $ij$. 
Together with the fact that the flip
graph is connected
(Lemma \ref{lem:props}(b))
and that $\phce_{f}$ has at least one vertex,
this implies that all pointed pseudo-triangulations appear
as vertices of $\phce_{f}$, and hence that all pointed and non-crossing
graphs appear as well.

Also, the extreme rays have the structure predicted in
Theorem~\ref{main}.
For a convex hull edge $ij$,
$K_{ij}$ must be an unbounded edge because
there is no other pointed pseudo-triangulation
that contains $E(v)-\{ij\}$.
\end{proof}

We now conclude that valid perturbation vectors $f\in
\reals^{\binom n 2}$ can be recognized by looking at
 $4$-point subsets only.

\begin{lemma}
\label{lem:4points}
A choice of $f\in \reals^{\binom n 2}$ is valid if and only if it is valid when
restricted to every four points of $P$.
\end{lemma}

\begin{proof}
By the previous Lemma, if $f$ is  not valid for $P$ 
then there is a point $v$ of  $\phce_{f}$
for which the graph
$E(v)$ is either non-pointed or crossing.
In either case,
there is a subset of four points $P'\subseteq P$
on which the induced subgraph is non-pointed or crossing.
Let $v'$ and $f'$ denote $v$ and $f$ restricted to $P'$.
Then, $v'$ is in $\phce_{f'}(P')$ and the graph $E(v')$ is crossing or not 
pointed, hence $f'$ is not valid on $P'$. Contradiction.
\end{proof}

\paragraph{\bf The case of four points.}
\begin{theorem}
\label{thm:R}
A choice of perturbation parameters $f\in \reals^{\binom 4 2}$ on four points
$P=(p_1,p_2,p_3,p_4)$ forms a valid choice if and only if 
\begin{equation}
\label{eq:sum_wf}
\sum_{1\leq i<j \leq 4} w_{ij} f_{ij} >0,
\end{equation}
where the $w_{ij}$'s are the unique self-stress on the four points, with
signs chosen as in Lemma~\ref{lem:stresses}.
\end{theorem}

For a set of four points $P=(p_1,p_2,p_3,p_4)$,
we denote by $G_{ij}$ the graph on $P$ whose only
missing edge is~$ij$. Recall that the choice of self-stress on four points
has the property
that $G_{ij}$ is pointed and non-crossing (equivalently, $ij$ is interior)
if and only if $w_{ij}$ is negative. 

Since $\phce_f(P)$ is five-dimensional, 
for every vertex $v$ the set $E(v)$ contains at least
five edges. Therefore $E(v)$ is either the complete graph or
one of the graphs $G_{ij}$. Theorem \ref{thm:R} is then a consequence of Lemma
\ref{lem:planar} and the following statement.

%\begin{lemma}
%\label{Lem:complete}
%The complete graph appears as a set of tight edges in $\pce_{f}$
%if and only if 
%\[
%%R=0. % 
%\sum_{1\leq i < j \leq 4} w_{ij} f_{ij} =0
%\]
%\end{lemma}
%
%\begin{proof} 
%Lemma~\ref{lem:stress4} gives the \emph{only if\/} part of 
%the statement.
%
%For the \emph{if\/} part, fix any particular edge $kl$ and consider
%a motion which produces as edge length increase $f_{ij}$ for every
%edge in $G_{kl}$. This
%motion exists (and is unique, modulo trivial motions) since every graph on 4
%non-degenerate points with 5 edges is infinitesimally minimally rigid.
%Lemma~\ref{lem:stress4} and
%our assumption  $\sum_{i,j} w_{ij} f_{ij}=0$
%imply that this motion gives the desired edge length increase for the edge
%$kl$ too. 
%\end{proof}

\begin{lemma}
\label{lem:incomplete}
Let $R := \sum_{1\leq i<j\leq 4} w_{ij} f_{ij}$.
For every edge $kl$, the following properties are equivalent:
\begin{enumerate}
\item The graph $G_{kl}$ appears as a vertex of $\pce_f(P)$.
\item $R$ and $w_{kl}$ have opposite signs.
\end{enumerate}
\end{lemma}

\begin{proof}
The graph $G_{kl}$ appears as a face 
if and only if the (unique, since $G_{kl}$ is rigid) 
motion with edge length increase $f_{ij}$
for every edge $ij$ other than $kl$ has edge length increase on $kl$ greater
than $f_{kl}$. In this case, by rigidity of $G_{kl}$, the face is actually a
vertex. But, for this motion:
\begin{eqnarray*}
0\!\!\! & = & \!\!\!\!\!
\sum_{1\leq i < j \leq 4} \!\!w_{ij} \langle p_j-p_i, v_j-v_i\rangle
= \! \sum_{1\leq i < j \leq 4} \!\! w_{ij} f_{ij} + 
w_{kl} (\langle p_k-p_l, v_k-v_l\rangle  - f_{kl})
\\
&
= & R + w_{kl} (\langle p_k-p_l, v_k-v_l\rangle - f_{kl}).
\end{eqnarray*}
Hence, $\langle p_k-p_l, v_k-v_l\rangle > f_{kl}$
is equivalent to $R$ and $w_{kl}$ having opposite sign.
\end{proof}

Observe that the previous lemma implicitly includes the statement that the
complete graph appears as a vertex if and only if $R=0$. 
The \emph{only if} part of this is actually an easy consequence of
Lemma~\ref{lem:stress4}. In this case $\pce_{f}$ degenerates to a single
point.

To complete the proof of Theorem \ref{main}
we still need to show that valid 
choices of perturbation parameters exist:

\begin{theorem}
\label{thm:valid-new}
Let $a$ and $b$ be any two points in the plane.
For any
point set $P=\{p_1,\ldots,p_n\}$ in general position in
the plane, the
following choice of parameters $f$ is valid:
\begin{equation}
f_{ij} = \det(a,p_i,p_j)\det(b,p_i,p_j)
\label{eq:ff'}
\end{equation}

\end{theorem}

\proof % \begin{proof}
This follows from the following Lemma, taking into account
  Lemma~\ref{lem:4points} and Theorem~\ref{thm:R}.

  \begin{lemma}
\label{lem:valid-2d}
For any two points $a$ and $b$ and
for any four points $p_1,\dots,p_4$ in general position in
the plane, we have
\[
\sum_{1\leq i<j\leq 4} w_{ij} f_{ij} % \det(a,p_i,p_j)\det(b,p_i,p_j)
 =1,
\]
where
$f_{ij}$ is given by~\eeqref{ff'} and
 $w_{ij}$ are the self-stress of Lemma~\ref{lem:stresses}.
  \end{lemma}
\proof % \begin{proof}
Let us consider
the four points $p_i$ as fixed and regard $R=\sum w_{ij} f_{ij}$ as a function
of $a$ and~$b$.
\[
R(a,b)=\sum_{1\leq i<j\leq 4} \det(a,p_i,p_j)\det(b,p_i,p_j) w_{ij}.
\]
For fixed $b$, $R(a,b)$ is clearly an affine function of $a$.
  We claim that $R(p_i,b)=1$
  for each of the four points $p_1,\dots,p_4$, which implies that $R(a,b)$ is
  constantly equal to $1$.
  To prove the claim, without loss of generality we take $a=p_1$.
Now, $R(p_1,b)$ is an affine function of $b$. By a similar argument
as before, it suffices to show
$R(p_1,b)=1$ for the three affinely independent points
$b=p_2,p_3,p_4$.
Without loss of generality we look only at $R(p_1,p_2)$.
Then, $f_{ij}=0$ for every $i,j$ except
$f_{34}=\det(p_1, p_3,p_4)\det(p_2, p_3,p_4)$.
Hence,
\[ 
R(p_1,p_2) = w_{34}f_{34} =
\frac
{\det(p_1, p_3,p_4)\det(p_2, p_3,p_4)}
{\det(p_3,p_4,p_1)\det(p_3,p_4,p_2)}
=1.
%\eqno 
\qed
\]
% %begin{eqnarray*}
%  &=&  \frac{\det(p_1, p_2,p_3)\det(b, p_2,p_3)}
%                  {\det(p_2,p_3,p_1) \det(p_2,p_3,p_4)}
% + \frac{\det(p_1, p_2,p_4)\det(b, p_2,p_4)}
%        {\det(p_2,p_4,p_1) \det(p_2,p_4,p_3)}
% + \frac{\det(p_1, p_3,p_4)\det(b, p_2,p_4)}
%        {\det(p_3,p_4,p_1) \det(p_3,p_4,p_2)} =
% \\
% & = 
% &\frac{\det(b, p_2,p_3)}{\det(p_2,p_3,p_4)}
% + \frac{\det(b, p_2,p_4)}{\det(p_2,p_4,p_3)}
% + \frac{\det(b, p_3,p_4)}{ \det(p_3,p_4,p_2)} =
% \\& =& \frac{1}{\det(p_2,p_3,p_4)}
%  ( \det(b, p_2,p_3)- \det(b, p_2,p_4) + \det(b, p_3,p_4) )= 1
% ,
% \end{eqnarray*}
% where the last equality follows from
% $\det(p_2,p_3,p_4)=
% \det(b, p_2,p_3) + \det(b, p_3,p_4) + \det(b, p_4,p_2).
% $
%\end{proof}

This proof is quite easy, but it does
not provide much intuition why this choice of $f$ works.
The first valid choice that we found by heuristic
considerations was the function
\begin{equation}
f'_{ij}={\textstyle \frac 12}\cdot
\left(|p_i|^2+|p_j|^2+\langle p_i,p_j\rangle\right)\cdot|p_i-p_j|^2
.
\label{eq:ff-old}
\end{equation}
The intuition behind this is as follows.
Looking back to Lemma~\ref{lem:4points}, there are two cases of four points:
in convex position and as a triangle with a point in the middle.
In both cases we want to avoid the situation that all interior
edges (inside the convex hull) are tight, while the hull edges
expand at least at their prescribed rate $f_{ij}$. Thus, we want
$f_{ij}$ to be big for the ``peripheral'' edges and
small for the ``central'' edges.
(This goal is in accordance with Theorem~\ref{thm:R}, as our
choice of  $\omega_{ij}$ is positive on boundary edges and negative on
interior ones.)

A function which has this property of being on average
bigger on the border of a region than in the middle is
the convex function $|x|^2$.
When we integrate $|x|^2$ over the edge $p_i p_j$
and multiply the result by the edge length (because
$f_{ij}$ is expressed in terms of the derivative of the \emph{squared}
edge length), we get~\eeqref{ff-old},
up to a multiplicative constant.

The parameters $f'_{ij}$ are valid.
Indeed, it can be checked that 
$\sum w_{ij} f'_{ij} = 1$ holds for all 4-tuples of points: setting
$a=b=0$ in the definition \eeqref{ff'} of $f$,
the difference
\[
f'_{ij}-f_{ij}=:
g_{ij}= 
\left[ |p_i|^2 |p_i|^2 + |p_j|^2 |p_j|^2 -  
        (|p_i|^2 + |p_j|^2)\langle p_i, p_j\rangle \right]/2
\]
satisfies
$\sum_{i,j} w_{ij} g_{ij} = 0$, which follows from
 Lemma \ref{lem:stress4} with $v_i= {|p_i|}^2 p_i/2$.

Of course, the equation
$\sum w_{ij} f'_{ij} = 1$
can be trivially
checked by expanding the values of $\omega_{ij}$ and $f'_{ij}$
with the help of a computer algebra package.
Attempts to find
a more classical proof have lead to
the function $f_{ij}$ of~\eeqref{ff'}.

%%%---------------------------

\section{Applications of the Main Result}
\seclab{appl}

\subsection{The Expansion Cone}
\seclab{x-cone}
As mentioned in the previous section,
the expansion cone is the recession cone $\phce_{0}$
of the pce-polyhedron $\phce_{f}$, whose structure we know.
The extreme rays of $\xcone$ are precisely
the extreme rays of  $\phce_{f}$, shifted to start at~0, but
parallel rays of  $\phce_{f}$ give rise to only
one ray of~$\xcone$, of course.

Studying when this happens will allow us to give now a rather easy proof of
the characterization of the extreme rays of~$\xcone$
(Proposition~\lemref{ex-rays}):
We conclude from Theorem~\ref{main} that the extreme rays correspond
to pointed pseudo-triangulations with one hull edge removed,
i.e., pte-mechanisms.
Any convex subpolygon in a pte-mechanism must be
rigid in the mechanism, according to
Lemma~\lemref{xcone-tight}b(iii).
This corresponds to the fact that every convex subpolygon
of a pointed pseudo-triangulation
contains
a pointed pseudo-triangulation
of that polygon and the enclosed points,
and is therefore rigid.
We still have to show that these r-components
are the only subcomponents that move rigidly
in the (unique) flex $v$ on a pte-mechanism~$G(P)$.

\begin{lemma}
 Let $P'\subset P$ be a maximal subset that moves 
 \textup(infinitesimally\textup) rigidly by
 the unique flex $v$ of the pte-mechanism $G(P)$.
 \begin{enumerate}
 \item [\rm(a)]
Then $P'$ contains all points of $P$ within its convex hull,
 \item [\rm(b)]
 $G$ contains no edge which either crosses the boundary
of the convex hull of $P'$ or gives a non-pointed graph together with the
boundary of $P'$, and
 \item [\rm(c)]
 $G$ contains all boundary edges of the
convex hull of~$P'$.
 \end{enumerate}
\end{lemma}
\begin{proof}
 (a)~A subset $P'\subset P$ moves rigidly
if and only if $E(v)$, considered in $\phce_{0}$,
 contains the complete subgraph
spanned by~$P'$. Then part~(a) follows from
Lemma~\lemref{xcone-tight}b(iii).

(b)~An edge $ij \in G\subseteq E(v)$ in either of those two situations
would, by Lemma~\lemref{xcone-tight}b(i) or~(ii), imply
that the complete graph on $P'\cup ij$
is part of $E(v)$. 
and hence $i$ and $j$ are rigidly connected to
$P'$. 
On the other hand, either of the two conditions implies
that one of $i$ or $j$ is not in $P'$, violating maximality of $P'$ as
a subset which moves rigidly.

(c)~Assume that a hull edge $ij$
of~$P'$ is missing in~$G$. Let $G'$ be the graph obtained adding $ij$ to $G$.
Since $P'$ moves rigidly, $G'$ is still a 1DOF mechanism. On the other hand, 
part (b) implies that $G'$
is still a pointed and non-crossing graph. Since $G$ had $2n-4$ edges,
$G'$ has $2n-3$ edges and, by Lemma \ref{lem:prop-pointed}, it is a pointed
pseudo-triangulation. This is a contradiction, because pointed
pseudo-triangulations are infinitesimally rigid
(Proposition~\ref{lem:pseudot}). 
\end{proof}

If follows from the last statement
that the rigidly moving components are precisely the convex
subpolygons of the pte-mechanism,
and two pte-mechanisms yield the same
motion (extreme ray) if and only if they
lead to the same collapsed pte-mechanism,
thus concluding the proof of
Proposition~\lemref{ex-rays}.
\qed

\subsection{Strictly Expansive Motions and Unfoldings of Polygons}

\begin{lemma} \lemlab{ex-motion}
Let 
$G(P)$ be a non-crossing and pointed framework in the plane.
Then, $G(P)$ has a non-trivial expansive flex if and only if it does not
contain all the convex hull edges. 
In this case, there is an expansive motion 
that is strictly expansive on all the convex hull edges not in $G$.
%$$\langle p_i-p_j, v_i-v_j\rangle > 0$$
\end{lemma}
\begin{proof}
  If all the convex hull edges are in $G$, then Lemma~\lemref{bounded} implies
  the statement: the face of $\phce_f$ corresponding to $G$ is bounded and,
  hence, it degenerates to the origin in $\phce_0$. If a certain convex hull
  edge $ij$ is not in $G$, then we extend $G$ to a pointed
  pseudo-triangulation, according to
  Lemma~\ref{lem:prop-pointed}. %\ref{maximal}). 
  Removing $ij$ yields a
  pte-mechanism, whose expansive motion is strictly expansive on $ij$. Adding
  all such motions for the different missing hull edges gives the stated
  motion.
\end{proof}

This immediately gives the following theorem.
\begin{theorem} \theolab{unfold-ex}
Let $G(P)$ be a non-crossing non-convex plane polygon
or a plane polygonal arc that does not lie on a straight line.
Then there is an expansive
motion that is strictly expansive on at least
one edge.
\qed
\end{theorem}

This statement has been crucial to show that every
simple polygon in the plane can be
unfolded by a global motion into convex position
and every polygonal arc can be straightened, without collisions
\cite{CDR,Streinu-2000}.
The proof given in those papers is based on
several reduction steps
between
infinitesimal motions, self-stresses, and polyhedral terrains.
The above new proof is completely independent, although not less
indirect.

Actually, one can work a little harder in the proof
of Lemma~\lemref{ex-motion} and show that
\emph{any} edge $ij\notin G$ that is not
contained inside a convex subpolygon of~$G$
can be chosen to be strictly expansive.
(The proof constructs an
appropriate pte-mechanism by a flipping argument, applied to the
minimal convex subpolygon enclosing the chosen edge.)
By adding several motions of $G(P)$ one can obtain
an expansive motion that is strictly expansive on \emph{all}
eligible edges, and hence, in Theorem~\theoref{unfold-ex},
there is a motion that is strictly expansive on all
edges $ij\notin G$.
This is actually the statement that was proved in~\cite{CDR},
in a more general setting.

\section{Other Constructions}
\label{other}

In this section we present three related results: a different
representation for the ppt-polytope that is less
dependent on some seemingly arbitrary choice of
parameters $f$, and the two constructions
leading to the associahedron: $2$-dimensional points in convex position
and the one-dimensional expansion polytope.

\subsection{A redefinition of the ppt-polyhedron}
\label{meta-ppt}
Let $P=\{p_1,\dots,p_n\}$ be a fixed point configuration
in general position.
As before, with each possible choice of parameters
$f=(f_{ij})\in \reals^{\binom n 2}$
we associate the polyhedron $\phce_{f}$ defined by the
constraints \eeqref{normal2} and
(\ref{eq:f}), and the polytope $\pce_{f}$ obtained setting to
equalities the inequalities corresponding to convex hull edges.
The case $f\equiv 0$ produces the expansion cone, with the polytope
degenerating to a point. The choice
$f_{ij}=\det(0,p_i,p_j)^2$, among others,
produces our polytope of
pointed pseudo-triangulations, according to Theorem \ref{thm:valid-new}.
But the results of Section \ref{polytope} imply
that actually any other choice of $f_{ij}$'s
would provide (combinatorially) the same polyhedron and polytope
as long as it satisfies inequality~\eeqref{sum_wf}
%\[
%  \sum_{i<j\in\{i_1,i_2,i_3,i_4\}} w_{ij} f_{ij} \ge 0
%\]
for every four points $p_{i_1}$, $p_{i_2}$, $p_{i_3}$ and $p_{i_4}$,
where the $w_{ij}$'s are the self-stress on the four
points with sign chosen as in Lemma \ref{lem:stresses}.

In this section we give a new construction for this
polyhedron, % of Section \ref{polytope}, 
with the advantage that it does
not depend on any choice of parameters. It has the disadvantage, however, that
it involves more variables: one for each of the ${\binom n 2}$ possible
edges among the $n$ points.

The basic idea is to study the image of the previously defined ppt-polytope
under the rigidity map 
 $M\colon \reals^{2n}\to  \reals^{\binom n 2}$
 for the complete graph on $P$,
using the variables
 $\delta_{ij} := \langle p_i-p_j, v_i-v_j\rangle$.
The following lemma is a stronger version of Lemma \ref{lem:stress4}.

\begin{lemma}
{\bf (Equivalence of parameters)}
\label{lem:equivalent-new}
Let $\delta=(\delta_{ij})$
be a vector in $\reals^{\binom n 2}$.
Then, the following properties are equivalent\textup:
\begin{enumerate}
\item[\textup{(a)}] $\delta\in \im M$
\item[\textup{(b)}]
For any four points $p_{i_1}$, $p_{i_2}$, $p_{i_3}$ and $p_{i_4}$
in $P$ one has
\[
  \sum_{i<j\in\{i_1,i_2,i_3,i_4\}} w_{ij} \delta_{ij} = 0
\]
where the $w_{ij}$'s are a nonzero self-stress of the complete graph on those
four points. 
\end{enumerate}
\end{lemma}

\begin{proof}

The implication (a)$\Rightarrow$(b) follows directly from one direction of
Lemma \ref{lem:stress4}.
Lemma \ref{lem:stress4} also gives (b)$\Rightarrow$(a) for each
quadruple: if (a) holds, then for each four
points there is an infinitesimal motion $(a_{i_1},a_{i_2},a_{i_3},a_{i_4})$
whose image by the rigidity map of the four points are the six relevant
entries of $\delta$. 
The motion for a quadruple is not unique, but any two
choices differ only by a trivial motion of the quadruple.

To define
a global motion $(a_1,\dots,a_n)$ of the whole configuration,
let us start by setting
$a_1=(0,0)$ and
$a_2=(0,b)$, where  $b$ must be the unique number satisfying
$\langle p_1-p_2, (0,b)\rangle = \delta_{12}$. (We assume without loss of
generality that $p_1$ and $p_2$ do not have the same $y$-coordinate.)
The condition $\delta=M(a)$ on the edges $1i$ and $2i$ then uniquely
defines $a_i$ for every $i\ne1,2$, because these two equations are
linearly independent, the directions $p_i-p_1$ and $p_i-p_2$
being not parallel.
To see that this global motion satisfies $\delta=M(a)$ also for any
other edge $ij$ ($i\ne1,2$, $j\ne1,2$), it is sufficient
to consider the quadruple $(p_1,p_2,p_i,p_j)$.
By construction,
the motion $(a_1,a_2,a_i,a_j)$ satisfies
the condition  $\delta=M(a)$ for five of the six edges in the quadruple.
Assumption~(b) says that $\sum_{k,l\in\{1,2,i,j\}}w_{kl}\delta_{ij}=0$
which, by Lemma \ref{lem:stress4}, implies that $\delta$ restricted
to $(p_1,p_2,p_i,p_j)$ is the set of edge increases produced by some motion
$v=(v_1,v_2,v_i,v_j)$. This motion can be normalized to $v_1=(0,0)$ and
$v_2=(0,b)$ and then it must coincide with $a$ by our uniqueness argument
above.
\end{proof}

This, together that the observation that the kernel of $M$ consists only of
 trivial motions  immediately gives the following lemma.
\begin{lemma}
\label{lem:meta-ppt}
For any $f\in \reals^{\binom n 2}$, the polyhedron $\phce_{f}$ is linearly
isomorphic to the one defined in ${\binom n 2}$
variables $\delta_{ij}=\delta_{ji}$
 by the following
${\binom n 4}$ equations and ${\binom n 2}$ inequalities.
\begin{eqnarray*}
&  \sum_{i<j\in\{i_1,i_2,i_3,i_4\}} w_{ij} \delta_{ij} = 0, \qquad
\forall\ i_1,i_2,i_3,i_4\in \{1,\dots,n\},
\\
& \delta_{ij}\ge f_{ij}, \qquad \forall\ i,j\in \{1,\dots,n\}.
\end{eqnarray*}
Moreover, setting to equalities the inequalities corresponding to convex hull
edges produces a polytope linearly isomorphic to $\pce_{f}$.
\qed
\end{lemma}

By making the change of variables $d_{ij}=f_{ij}-\delta_{ij}$  and
taking into account that $\sum w_{ij} f_{ij}=1$ for any four points, for the
valid choices $f$ of Theorem \ref{thm:valid-new},
we conclude:
\begin{theorem}
\label{theo:meta-ppt}
For any given point set $P=\{p_1,\dots,p_n\}$ in the plane in general
position, the
following ${\binom n 4}$ equations and ${\binom n 2}$
inequalities define a
simple polyhedron in $\reals^{\binom n2}$
linearly isomorphic to the polyhedron $\phce_f(P)$ of Theorem~\ref{main},
with $f_{ij}$ chosen as in Theorem~\ref{thm:valid-new}.
\begin{itemize}
\item %[\rm (a)] 
$\sum w_{ij} d_{ij} =1$ for every quadruple,
where the $w_{ij}$'s of each equation are
the self-stress on the corresponding quadruple
stated in Lemma \ref{lem:stresses}.
\item %[\rm (b)] 
$d_{ij}\le 0$ for every variable.
\end{itemize}
The maximal bounded face in the polyhedron
is obtained setting to equalities the inequalities corresponding
to convex hull edges.
\qed
\end{theorem}
The $\binom n 4$ equations % in part (a) 
are highly redundant.
It follows from the proof of Lemma~\ref{lem:equivalent-new}
that the $\binom{n-2} 2$ quadruples involving two fixed vertices
are sufficient. Subtracting this number from the number
$\binom n 2$ of variables actually gives the right dimension $2n-3$
of the polyhedron.

\subsection{Convex position and the associahedron}
\label{convex2dim}

Suppose now that our $n$ points $P=\{p_1,\dots,p_n\}$ are in (ordered) convex
position.
In this subsection all indices are regarded modulo~$n$.
In Section \ref{preliminaries} we noticed that the polytope of pointed
pseudo-triangulations of $P$ is combinatorially the same thing as the
secondary polytope, which for a convex $n$-gon is an associahedron. We prove
here that, in fact, the secondary polytope and the ppt-polytope
are affinely isomorphic.

The first problem we encounter is that so far we have only facet descriptions
for the ppt-polytope, while the secondary polytope is defined by the
coordinates of its vertices. We recall that the secondary polytope lives in
$\reals^n$ and that the $i$-th coordinate of the vertex corresponding to a
certain triangulation $T$ equals the total area of all triangles of $T$
incident to $p_i$. Denote this area as $\Area_T(p_i)$. For convenience
we will work with a normalized definition of area of a triangle with vertices
$p$, $q$ and $r$ as being equal to $|\det(p,q,r)|$.
We also assume our points to be ordered
counter-clockwise, so that $\det(p_i,p_j,p_k)$ is positive if and only if $i$, $j$ and
$k$ appear in this order in the cyclic ordering of $\{1,\dots,n\}$.
In this way
\[
\Area_T(p_i) := \sum_{l=1}^{t-1} \det(p_i,p_{j_l},p_{j_{l+1}})
\]
where $\{p_{j_1},\dots,p_{j_t}\}$ is the ordered sequence of vertices adjacent
to $p_i$ in $T$.

Our first task is to compute the coordinates of the vertices of the
ppt-polytope. Notice that by definition the coordinates corresponding to edges
of $T$ are zero, since the inequalities corresponding to the edges of $T$ are
satisfied with equality. It will turn out that we do not need all the
coordinates, but only those corresponding to \emph{almost-convex-hull edges}
$p_{i-1}p_{i+1}$.

\begin{lemma}
\label{lemma:convex}
Let $T$ be a triangulation of $P$. Then, in the ppt-polytope of Theorem
\ref{theo:meta-ppt}, the coordinate $d_{i-1,i+1}$ corresponding to an
almost-convex-hull edge equals
\[
d_{i-1,i+1} = - \det(p_{i-1},p_i,p_{i+1})
          \left(\Area_T(p_i)- \det(p_{i-1},p_i,p_{i+1}) \right).
\]
\end{lemma}

\begin{proof}
Let $p_{j_1},\dots,p_{j_t}$ be the ordered sequence of vertices adjacent to
$p_i$ in $T$, with $p_{j_1}=p_{i+1}$ and $p_{j_t}=p_{i-1}$. We will
prove by induction on $k=2,\ldots,t$ that
the coordinate
$d_{j_1 j_k}$ of the ppt-polytope vertex corresponding to $T$ equals
\[
 d_{j_1 j_k} = - \det(p_i,p_{j_1},p_{j_k}) \Area(p_{j_1},\dots,p_{j_k})
\]
where $\Area (p_{j_1},\dots,p_{j_k})$ denotes the area of the polygon with
vertices $p_{j_1},\dots,\allowbreak p_{j_k}$, in that order. This
reduces to the formula in the statement for $k=t$.

The base case $k=2$ is trivial: since $j_1 j_2$ is an edge in
the triangulation, we have $d_{j_1 j_2}=0$.
To compute $d_{j_1j_{k}}$ for $k>2$ we consider the
quadruple $p_i, p_{j_1}, p_{j_{k-1}}, p_{j_{k}}$.
The only non-zero $d$'s on this quadruple are $d_{j_1 j_{k-1}}$ and
$d_{j_{1} j_{k}}$.
(For $k=3$,  $d_{j_1 j_{k-1}}$ is also zero.)
Hence, the equation
$\sum w_{\alpha\beta} d_{\alpha\beta}=1$
for this quadruple reduces to
\[
d_{j_1 j_k} w_{j_1 j_k} + d_{j_{1} j_{k-1}} w_{j_1 j_{k-1}} =1.
\]
 From this we infer the stated value for $d_{j_1 j_k}$ from the known
values of the other quantities:
\begin{eqnarray*}
   w_{j_1 j_k}&=& \frac{1}{\det(p_i,p_{j_1},p_{j_k})
                \det(p_{j_1},p_{j_{k-1}},p_{j_k})},\\
   w_{j_{1} j_{k-1}}&=& \frac{-1}{\det(p_i,p_{j_{1}},p_{j_{k-1}})
                \det(p_{j_1},p_{j_{k-1}},p_{j_k})},
 \hbox{\ and}\\
   d_{j_{1} j_{k-1}} &=& -\det(p_i,p_{j_{1}},p_{j_{k-1}})
                  \Area(p_{j_{1}},\dots,p_{j_{k-1}}).
\end{eqnarray*}
(The last equation is the inductive hypothesis.)
\end{proof}

This Lemma immediately implies the following:

\begin{corollary}
\label{coro:convex}
\begin{enumerate}
\item[\rm (a)]% \hskip7pt] \null\hskip-9pt 
The affine map
$$
a_i = -
   \frac { d_{i-1,i+1}}{ \det(p_{i-1},p_{i},p_{i+1})} +
   \det(p_{i-1},p_{i},p_{i+1})
$$
gives the
coordinates $(a_1,\dots,a_n)$
of a triangulation %~$T$
in the secondary polytope 
 in terms of its coordinates
$(d_{ij})_{ij\in{\binom n 2}}$ % of $T$
 in the ppt-polytope of
Theorem~\ref{theo:meta-ppt}.
\item[\rm (b)] The substitution $d_{ij}=f_{ij} - \langle p_i-p_j, v_i-v_j\rangle$
in the above formula gives an affine map $\reals^{2n}\to \reals^n$ sending the
ppt-polytope of Theorem~\ref{main} to the secondary polytope, whenever the
$f_{ij}$'s are a valid choice satisfying 
\[
\sum_{1\leq i < j \leq 4} w_{ij} f_{ij}=1
\]
for every quadruple, such as the choices in Theorem \ref{thm:valid-new}.
\qed
\end{enumerate}
\end{corollary}

Observe that Corollary \ref{coro:convex} implies that, for points in convex
position, we can consider the ppt-polytope of Theorem \ref{theo:meta-ppt}
as lying in the $n$-dimensional space given by the
coordinates $d_{i-1,i+1}$.
For the polyhedron, we
additionally need the coordinates $d_{i,i+1}$. (These are zero on the polytope,
but not on the polyhedron.)

%
%Remark: The corollaries imply that we can consider the mpt-metapolytope as
%lying in the n-dimensional projection of $\reals^{n\choose 2}$ given by
%the coordinates $d_{i-1,i+1}$. This is not true if we want the
%mpt-metapolyhedron (which has dimension $2n-3$). For this we additionally
%need the n coordinates d_{i,i+1}, which are identically zero in the
%polytope but not in the polyhedron.
%
%Actually, it is obvious in general (i.e., even for points not in
%convex position) that knowing the $2n-3$ $d_{ij}$'s corresponding to any
%spanning and rigid graph fixes all the other $d_{ij}$'s (although the
%explicit formulas for them are far from trivial).%
%
%Guenter, is it easy for you to check:
%(1) whether the formulas given in Corollary 2 are the same as yours or
%(2) they at least produce a linear isomorphism as claimed?%
%

%\section{Towards a general theory in arbitrary dimension}
%\label{arbitrary_dim}

%We can
%express this and the fact that the kernel of the rigidity map equals
%the trivial motions (motions coming
%from rigid motions of the whole plane) by saying that the
%following sequence of linear maps is \emph{exact},
%\compaint{these sequence needs better typesetting. Paco}
%\[
%0 \to \reals^3 \to^i\to \reals^{2n} \to^M\to \reals^{n\choose 2}
%\to^\Delta\to \reals^{n\choose 4}
%\]
%where $i$ is just the inclusion of the linear space of trivial motions (which
%is 3-dimensional) in the space of all motions.

\subsection{The 1-Dimensional Case: the Associahedron Again}
\label{1dim}

By considering $1$-dimensional
expansive motions, in this section we will recover the associahedron via a
different route. The analogy of this construction to the 2-dimensional case
will become even more apparent in Section \ref{towards}.

\paragraph{\bf The polytope of constrained expansions in dimension 1.}
\label{sec:exp}
In the 1-dimensional case we will rewrite
constraints (\ref{eq:f}) in the form
\begin{equation}
\label{eq:f1}
v_j - v_i \ge g_{ij},\qquad \forall i < j.
\end{equation}
One set of inequalities is equivalent to the other under the change of
constants $g_{ij}(p_j-p_i) = f_{ij}$, for any $i<j$. This reformulation
explicitly shows that the solution set does not depend on the point set
$P=\{p_1,\dots,p_n\}$ 
that we choose. We denote this solution set $\phce_{g}$,
to mimic the notation of the 2-dimensional case.

It is easy to see that the polyhedron $\phce_{g}$ 
is full-dimensional and 
it contains no lines if we add the normalization equation $v_1=0$.
Hence, after normalization, it has dimension $n-1$ and contains some vertex.
For any vertex $v$ or for any feasible point $v\in \phce_{g}$,
we may look at the set $E(v)$ of tight inequalities at~$v$:
$$
 E(v) := \{\,ij \mid 1\le i<j\le n,  v_j-v_i = g_{ij}\,\}
$$
We regard $E(v)$ as the set of edges of a graph on the
vertices $\{1,\ldots,n\}$.

One may get various polyhedra by choosing different numbers
$g_{ij}$ in~(\ref{eq:f1}).
We choose them with the following properties.
\begin{equation} \label{eq:f4}
g_{il}+g_{jk} > g_{ik} + g_{jl} ,
\qquad\forall 1\le i< j\le k<l\le n.
\end{equation}
For $j=k$ we use this with the interpretation $g_{jj}=0$, so
we require
\begin{equation} \label{eq:f3}
g_{il}  > g_{ik} + g_{kl} ,
\qquad\forall 1\le i< k<l\le n.
\end{equation}
One way to satisfy these conditions is to select
\begin{equation}
g_{ij} := h(t_j-t_i), \quad\forall i<j \label{eq:h}
\end{equation}
for an arbitrary strictly convex function $h$ with $h(0)=0$
and arbitrary real numbers $t_1 < \cdots < t_n$.
The simplest choice is $h(t)=t^2$ and $t_i=i$, yielding
$g_{ij} = (i-j)^2$.

Two edges $ij$ and $jk$ with $i<j<k$ are called
 \emph{transitive edges},
and edges $ik$ and $jl$ with $i<j<k<l$ are called
 \emph{crossing edges}.

\begin{lemma} \label{lem:valid}
If $g$ satisfies $(\ref{eq:f4}$--$\ref{eq:f3})$ and $v\in \phce_{g}$, then
$E(v)$ cannot contain {transitive}
or {crossing} edges.
\end{lemma}
\begin{proof}
If we have two transitive edges $ij,jk\in E(v)$ this means that
$ v_j-v_i = g_{ij}$ and
$ v_k-v_j = g_{jk}$.
This gives
$ v_k-v_i = g_{ij}+ g_{jk} < g_{ik}$,
by~(\ref{eq:f3}), and thus $v$ cannot be in~$\phce_{g}$
because it violates~(\ref{eq:f}).
The other statement follows similarly.
\end{proof}

\paragraph{\bf Non-crossing alternating trees.}

A graph without transitive edges is called an \emph{alternating}
or \emph{intransitive} graph:
every path in an alternating graph changes continually between up and down.

\begin{lemma}
A graph on the vertex set $\{1,\ldots,n\}$
 without transitive or crossing edges cannot contain
a cycle.
\end{lemma}
\begin{proof}
Assume that $C$ is a cycle without transitive edges.
Let $i$ and $m$ be the lowest and the highest-numbered
vertex of a cycle $C$, and let $ik$ be an edge of $C$ incident to $i$,
but different from $im$. The next vertex on the cycle after $k$ must
be between $i$ and $k$;
continuing the cycle, we must eventually reach
$m$, so there must be an edge $jl$ which jumps over $k$,
and we have a pair $ik$, $jl$ of crossing edges.
\end{proof}

Since the polyhedron is $(n-1)$-dimensional, the set $E(v)$ of a vertex
$v$ must contain at least $n-1$ edges. We have just seen that it
is acyclic, and hence it must be a tree and contain exactly
$n-1$ edges. So we get:
\begin{proposition} \label{prop:necessary}
$\phce_{g}$ is a simple polyhedron.
The tight inequalities for each vertex correspond to
non-crossing alternating trees.
\qed
\end{proposition}
We will see below that $\phce_{g}$ contains in fact \emph{all}
non-crossing alternating trees as vertices.

\paragraph{\bf A new realization of the associahedron.}
Let's look at the combinatorial properties of
these trees.
Alternating trees have been studied in combinatorics
in several papers,
see for example~\cite{Post,PS} or~
\cite[Exercise 5.41, pp.~90--92]{Stanley2}
and the references
given there.
Non-crossing alternating trees were studied only
by Gelfand, Graev, and Postnikov, %~\cite[Section 6]{GGP},
under the name of ``standard trees''. They proved the following
fact~\cite[Theorem~6.4]{GGP}.

\begin{proposition}
The non-crossing alternating trees non $n+1$ points are in one-to-one
correspondence with the binary trees on $n$ vertices,
and hence their number is the $n$-th Catalan number
${\binom{2n} n}/(n+1)$.
\qed
\end{proposition}

The bijection given in~\cite{GGP} to prove this fact
is that the vertices of the binary tree correspond to the edges of
the alternating tree. It is easy to see that every
non-crossing alternating tree must contain the edge $1n$.
Removing this edge splits the tree into two parts;
they correspond to the two subtrees of the root in the binary tree.
The two parts are handled recursively.
Even simpler is the bijection to bracketings (ways to insert $n-1$ pairs of 
parentheses in a string of $n$ letters).
Just change edge $ij$
by a parenthesis enclosing the $i$-th and $j$-th letter.

\begin{figure}[htb]
  \centering

\medskip

   \includegraphics[width=\textwidth]{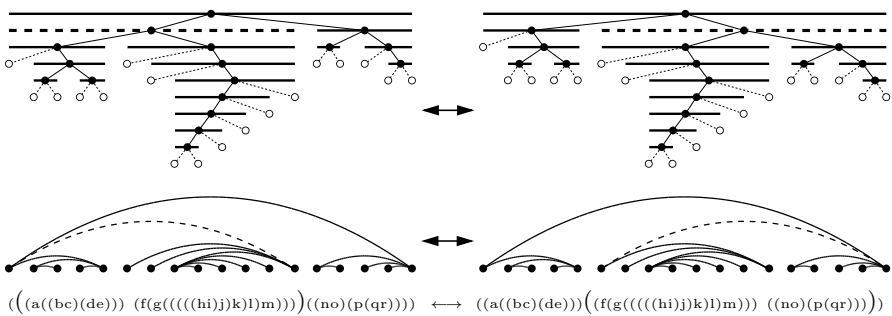}
%  \Ipe{bas2.ipe}

\medskip

\tiny
({\small(}(a((bc)(de))) (f(g(((((hi)j)k)l)m))){\small)}((no)(p(qr))))
$\ \longleftrightarrow\ $
((a((bc)(de))){\small(}(f(g(((((hi)j)k)l)m))) ((no)(p(qr))){\small)})
    \caption{The bijection between binary trees (up), non-crossing
alternating trees (middle) and bracketings (bottom). The flipping
 operation in each case is shown, and the elements affected by the
flip are highlighted.}
\label{fig:rotate}
\end{figure}

Figure~\ref{fig:rotate} gives an example of these correspondences, including
the correspondence of the respective flipping operations:
%
%We extend this correspondence to the adjacency structure between trees:
%\begin{lemma} \label{lem:flip}
If we remove any edge $e\ne 1n$ from a
 non-crossing alternating tree
$T$, there is precisely one other
 non-crossing alternating tree
$T'$ which shares the edges $T-\{e\}$ with~$T$.
This  exchange operation 
corresponds  under the above bijection
to a tree rotation of the binary tree, and to a single application of the
associative law (remove a pair of interior 
parentheses  and insert another one in the only possible way) in a
bracketing.

Observe that still we don't know that {\em all} the non-crossing and
alternating trees appear as vertices of $\pce_g$. But this is easy to prove:
since $\pce_g$ is simple of dimension $d-2$, its graph is regular of degree
$d-2$. And it is a subgraph of the graph of rotations between binary trees,
which is also $(d-2)$-regular and connected. Hence, the two graphs coincide.

%\qed
%\end{lemma}

If we now look at faces of $\phce_g$, rather than vertices, the tight edges
for each of them form a non-crossing alternating {\em forest}.
Such a forest $G$ is an expansive mechanism if and only if the edge $1n$ is
not in $G$: If that edge is present, $p_1$ and $p_n$ are fixed and any other
$p_i$ must approach one of the two. If the edge $1n$ is not present, then let
$i$ be the maximum index for which $1i$ is present. Then the motion 
$v_1=\cdots=v_i=0<v_{i+1}=\cdots=v_n$ is an expansive flex of $G$.

Hence, $\phce_g$ has a unique maximal bounded face,
the facet given by the equation
%\begin{equation} \label{eq:long}
$  v_n-v_1 = g_{1n}$ 
%\end{equation}
and corresponding to the edge $1n$ alone. This facet is then an
$(n-2)$-dimensional simple polytope that we denote $\pce_g$.
The $n-2$ neighbors of each vertex $v$ correspond to the $n-2$ possible
exchanges of edges different from $1n$ in the tree $E(v)$.

\begin{theorem} \label{thm:dim1}
$\pce_{g}$ is a simple $(n-2)$-dimensional 
polytope whose face poset is that of the associahedron.
Vertices are in one-to-one
correspondence with the non-crossing alternating trees
on $n$ vertices.
%or with the binary trees on $n-1$ vertices.
Two vertices are adjacent if and only if the two
 non-crossing alternating trees
differ in a single edge.
%$($or if the two binary trees differ by a rotation$)$.
%Hence it is an associahedron.

$\phce_{g}$ is an unbounded  polyhedron with the same vertex set 
as~$\pce_{g}$.
The extreme rays correspond to the
 non-crossing alternating trees
with the edge $1n$ removed.
\end{theorem}

\begin{proof}
  Only the statement regarding the face poset remains to be proved. This can
  be proved in two ways: On the one hand, we already know that the graphs of
  $\pce_g$ and of the $(d-2)$-associahedron coincide (the latter being the
  graph of rotations between binary trees). And simple polytopes with the same
  graph have also the same face poset. (This is a result of Blind and Mani; see
  \cite[Section 3.4]{Ziegler}). As a second proof, the correspondence between
  non-crossing alternating trees and bracketings trivially extends to a
  correspondence between non-crossing alternating forests containing $1n$ and
  ``partial bracketings'' in a string of $n$ letters which include the
  parentheses enclosing the whole string. The poset of such things is the face
  poset of the associahedron (see~\cite{Stanley2}, Proposition 6.2.1 and
  Exercise~6.33). In particular, the face-poset of $\pce_g$ is a subposet of
  the face-poset of the associahedron, and two polytopes of the same dimension
  cannot have their face posets properly contained in one another. This is 
  true in general by  topological reasons, but specially obvious in our case
  since we know our polytopes to be
  simple and their vertices to correspond one to one. Each vertex is incident
  to exactly $\binom{d-2}i$ faces of dimension $i$ in both polytopes, and
  two subsets of cardinality $\binom{d-2} i$ cannot be properly contained in
  one another.
%
%We prove only the statements regarding~$\pce_{g}$.
%We know that it has at least one vertex.
%By Proposition~\ref{prop:necessary}, that vertex
%must correspond to a
% non-crossing alternating tree.
%By Lemma~\ref{lem:connected}, every
%exchange neighbor of a
% non-crossing alternating tree
%that is represented in the polytope
%must also appear as a vertex.
%Since the set of all
% non-crossing alternating trees
%is connected under the edge exchange operation
%(like the set of binary trees under tree rotations),
%we conclude that all trees appear as vertices.
\end{proof}

%We remark that we have obtained this result in a somewhat indirect way,
%by combining combinatorial properties
%with general structural knowledge about simple polytopes.
%We have not explicitly proved that any single tree $E(v)$ is in fact
%feasible, i.~e., satisfies the constraints~(\ref{eq:f}).

A result which is related to Theorem~\ref{thm:dim1} was proved
by Gelfand, Graev, and Postnikov~\cite[Theorem~6.3]{GGP},
in a setting dual to ours:
there a triangulation of a certain polytope was constructed. The
 non-crossing alternating trees
correspond to the \emph{simplices}
of the triangulation.
It is shown explicitly that the simplices form a partition of the
polytope.
Certain numbers $g_{ij}$ are then associated with the \emph{vertices} of
the polytope to show that the triangulation is a projection of
the boundary of a higher-dimensional polytope.
Incidentally, the numbers that were suggested for this purpose are
$(i-j)^2$, which coincides with the simple proposal given above
after~\eeqref{h}, but the calculations are not given
in the paper~\cite{GGP}. 

One easily sees that
the conditions $(\ref{eq:f4}$--$\ref{eq:f3})$ on $g$ are also
necessary for the theorem to hold:
If any of these conditions would hold as an equality or
as an inequality in the opposite direction, the argument
of Lemma~\ref{lem:valid} would work in the opposite direction,
and certain
 non-crossing alternating trees
would be excluded.
Thus, $(\ref{eq:f4}$--$\ref{eq:f3})$ are a complete characterization
of the ``valid'' parameter values~$g_{ij}$.

\paragraph{\bf Further Remarks.}
The result we presented in this section is surprising in two ways:
first, that it produces such a well studied
object as the associahedron; second, that it requires additional types of
linear constraints that are not needed in dimension $2$. Indeed, inequalities
\eeqref{f3} in the 1-dimensional case are the exact analog of inequalities
\eeqref{sum_wf} in the 2-dimensional case, as we will see in Section
\ref{towards}. But the constraints \eeqref{f4} have no
analog. This second aspect
makes the task of producing $3$d generalization of the constructions of this
paper more
challenging, as there does not seem to be a straightforward pattern for
producing the linear constraints whose
instantiations in 1d and 2d give the polytopes of expansive motions.

The conditions (\ref{eq:f4}--\ref{eq:f3}) leave a lot of freedom for the
choice of the variables~$g_{ij}$.
We have an ${\binom n 2}$-dimensional
parameter space.  This is in contrast to the classical representation,
which has $2n$ parameters (the coordinates of
$n$ points in the plane).
If we select the parameters to be integral, we obtain an associahedron which
is an integral polytope (this is also true for the classical
associahedron for a polygon with integer vertices.)
But observe that, in fact, the associahedra
obtained here are in a sense much more special than the classical associahedra
obtained as secondary polytopes.  They have $n-2$ pairs of parallel facets,
given by the equations $v_i=v_1+g_{1i}$ and $v_i=v_n-g_{in}$ (i.e.,
corresponding to the pairs of edges $1i$ and $in$). 
See the right picture of Figure~\ref{fig:assoc}, where $n=5$.
This is no surprise, since
we have constructed our polytope by perturbing (a region of) a Coxeter
arrangement, whose hyperplanes are in extremely non-general position.

One can check that
classical associahedra have no parallel facets. Hence,
are associahedra are not affinely equivalent to the classical ones.
Recently, Chapoton, Fomin and Zelevinsky 
\cite{ChFoZe} have shown another construction of
the associahedron having as facet normals the root
system of type $A_n$ (more generally, they show a similar
construction for all root systems). Still, our associahedra
are not affinely equivalent to theirs: Comparing the right
part of our Figure~\ref{fig:assoc} with their Figure 2
we see that both are 3-dimensional associahedra with 
three pairs of parallel facets but in their realization the
three remaining facets share a vertex while in ours they do not.
It is however conceivable that they are related
by projective transformations.

\begin{figure}[htb]
  \centering
{\def\IPEfile{2d-assocs.ipe}\input{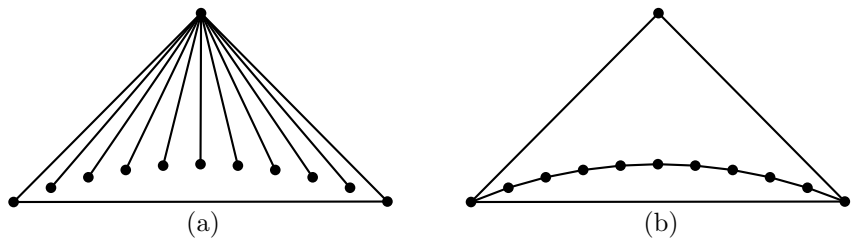}}
\caption{
(a)~The non-crossing alternating trees appear as
a face of a ppt-polytope.
(b)~Another associahedron in the same ppt-polytope.
}
\figlab{2d-assocs}
\end{figure}

The non-crossing alternating trees also appear directly
as a face of the ppt-polytope. Consider the configuration
of points shown in Figure~\figref{2d-assocs}, consisting
of a convex chain $\{p_1,\ldots,p_n\}$ and an additional point~$p_0$.
The pseudotriangulations containing
all the edges from $p_0$ to the other points
(Figure~\figref{2d-assocs}a), correspond to the
 non-crossing alternating trees on $p_1,\ldots,p_n$,
forming an associahedron of dimension~$n-2$.
Another associahedron (of dimension $n-3$)
 appears as a different face
of this ppt-polytope: the face corresponding to the edges
of the convex chain (Figure~\figref{2d-assocs}b).
These two faces are disjoint.

\paragraph{\bf Relations to optimization and the Monge matrices.}
Constraints of
the form~\eeqref{f1} are commonly found
in project planning and critical path analysis.
The variables $v_i$ represent unknown start times of tasks,
and the constraints~\eeqref{f1} specify waiting conditions
between different tasks.
For example, the minimization of $v_n-v_1$
% subject to the conditions
is a longest path problem in an acyclic network.
The different vertices of $\pce_g$ correspond to the
different optimal solutions when we apply various
linear objective functions.

Matrices $G=(g_{ij})$ with the property~\eeqref{f4} are said
to have the \emph{Monge property}, if we set $g_{ii}=0$
and $g_{ij}=\infty$ for $i>j$. The Monge property has received
a great deal of attention in optimization because it arises
often in applications and it characterizes special classes of
optimization problems that can be solved efficiently,
see~\cite{Monge} for a survey.

The dual linear program of the linear programming problem~\eeqref{f1}
(with a suitable objective function) is a minimum-cost flow problem
on an acyclic network with edges $ij$ for $1\le i<j\le n$, and
cost coefficients~$-g_{ij}$. Network flow is actually one of
the oldest areas in optimization in which the Monge property has been applied,
and where it has been shown that optimal solutions can be obtained
by a greedy algorithm in certain cases.
The non-crossing alternating trees are just the different possible subgraphs
of those edges which carry flow in an optimal solution.

\section{Towards a General Framework, in Arbitrary Dimension}
\label{towards}
To make more explicit the analogy between our 1d and 2d
constructions we'll rewrite the $1$d construction back in terms of
inequalities \eeqref{f} instead of \eeqref{f1}. As we mentioned 
at the beginning of Section \ref{1dim}, the way to
do this is to substitute $g_{i,j}=f_{i,j}/(p_j-p_i)$ everywhere. In particular,
the constraints \eeqref{f3} become
\[
\frac{f_{ik} }{ p_i-p_k} + \frac{f_{kl} }{ p_k-p_l} + \frac{f_{il} }{ p_l-p_i} >0.
\]
But $w_{ik}=\frac{1 }{ p_i-p_k}$, 
$w_{kl}=\frac{1 }{ p_k-p_l}$ and $w_{il}=\frac{1}{p_l-p_i}$ define a
self-stress 
on any 1-dimensional point set $\{p_i,p_k,p_l\}$. Hence, 
inequalities \eeqref{f3} are the exact analog
of inequalities (\ref{eq:sum_wf}) of the 2d case. The difference is that in 2d
these inequalities are necessary and sufficient for a choice of parameters to be
``valid'', while in 1d we need the additional constraints (\ref{eq:f4}), which do
not follow from (\ref{eq:f3}) as the following example shows:
\[
g_{12}=g_{23}=g_{34}=1,\qquad
g_{13}=g_{24}=2.2,\qquad
g_{14}=3.3.
\]

Let's see what would be a generalization of the 1d and 2d cases to arbitrary
dimension. 
For any finite point set $P=\{p_1,\dots,p_n\}$ spanning
$\reals^d$ (in general position or not) there is a well-defined cone
of expansive motions $\xcone$ of dimension $nd-\binom{d+1}2$, 
which is the positive region in the arrangement of
hyperplanes defined by the constraints~(\ref{eq:0}). Every choice of perturbation
parameters $(f_{ij})_{i,j\in\{1,\dots,n\}}$, via inequalities~(\ref{eq:f}),
produces a polyhedron of constrained expansions $\phce_f$,
whose faces are in
correspondence with some graphs embedded on $P$. If the $f_{ij}$'s are
sufficiently generic, this polyhedron is simple and:
\begin{itemize}
\item The face poset of $\phce_f$ is isomorphic to a 
  subposet of the graphs
  embedded on $P$, ordered by reverse inclusion.
\item The maximal graphs in this subposet are some of the infinitesimally
  minimally rigid frameworks on~$P$.
\item Adjacent vertices on the polyhedron correspond to
graphs differing only by one edge.
\end{itemize}

The main question is whether there is a clever choice of the $f_{ij}$'s or a
sensible choice of a special
family of minimally rigid frameworks to be the vertices
of the polyhedron which admits a nice and simple geometric characterization, 
as non-crossing alternating trees in 1d or pointed pseudo-triangulations in 2d
do. In particular, what should be the analog of pointed
pseudo-triangulations on a point set in general position in 3d? 
One difficulty is that the combinatorics of minimally rigid
graphs in dimension three is not fully understood.  Perhaps a good case to
start with would be points in convex and general position in dimension $3$, for
which the convex hull edges provide a very canonical minimally rigid
framework to start with.

Another observation is that both in the 1d and 2d cases
that we have treated here, the resulting polytope $\phce_f$ has a unique
maximal bounded face. This is related to the graph of flips between the
special minimally rigid frameworks considered being connected and
regular,
which is an extra desirable property.

It would even be interesting to have a good generalization of the 2d case to
point sets which are not in general position.
If we use the 2d definition of $f$ from equations \eeqref{ff'}
with many points on a common line,
there are solutions in which all the inequalities are tight.
In this way we get essentially the one-dimensional expansion cone
of Proposition~\lemref{X1}, when we project all vectors
$v_i$ on the direction of that line.
(This is the same situation as when all
relations~\eeqref{f3} are satisfied as equations.)
One way to get rid of this degeneracy is to
``perturb the perturbations'' by adding an infinitesimal
component of the one-dimensional expansion parameters 
of Section~\ref{1dim}:
instead of 
$f$ we use $f^{(2)}_{ij} + \eps f^{(1)}_{ij}$,
where  $f^{(2)}$ is valid for two-dimensional point sets, % given by~\eeqref{ff},
$f^{(1)}_{ij}$ satisfy the restrictions required for the 1d perturbation
parameters, and $\eps>0$ is sufficiently small.
The situation would be like in
Figure~\figref{2d-assocs}a when the points on the convex
chain converge to being collinear.
We have not further investigated this idea.

%\bigskip

Another consideration which may help to treat degenerate
point sets and arbitrary-dimensional ones is that the ideas in Section
\ref{meta-ppt} generalize with not much difficulty under very weak general
position assumptions. Let us first try to generalize 
Lemma~\ref{lem:equivalent-new}. 
We need to decide what is the right analog
of the statement in part (b). We propose the following.
Remember that a \emph{circuit} of a point set 
$P=\{p_1,\dots,p_n\}\subset\reals^d$ is a
minimal affinely dependent subset. If $P$ is in general position, its set of
circuits ${\cal C}(P)$ is $\binom P{d+2}$.
In general, any set of $d+2$ points spanning $d$ dimensions contains a unique
circuit, but this circuit can have less than $d+2$ elements. 
In any case, the complete graph on a
% circuit has the
% property that the removal of any single edge provides a minimally rigid
% framework. In particular, each
circuit $C$ has (up to a constant factor) 
a unique self stress $(w_{ij})_{i,j\in C}$ with no vanishing $w_{ij}$, which
is the one given by Lemma \ref{lem:dependent}.

Let us consider the following linear maps:
\begin{eqnarray*}
M\colon \qquad \reals^{dn}\quad &\longrightarrow & \quad\reals^{\binom n 2}\\
     (v_1,\dots,v_n) &   \longmapsto   &
                          \langle v_i-v_j, p_i -p_j \rangle,
\\[1ex]
%\end{eqnarray*}
%\begin{eqnarray*}
\Delta\colon \qquad\reals^{\binom n 2}\quad
                   &\longrightarrow&
                           \qquad\reals^{{\cal C}(P)}\\
  \{f_{ij}\}_{ij\in{\binom n 2}}
                    &   \longmapsto   &
                \Big\{\sum_{ij\in C} w_{ij}f_{ij}\Big\}_{C\in{\cal C}(P)}.
\end{eqnarray*}
$M$ is just the rigidity map of the complete graph on $P$. The equivalence of
(b) and (c) in Lemma~\ref{lem:equivalent-new} can be rephrased as
\[
\ker \Delta=\im M, \hbox{ for every planar point set in general position.}
\]
But in fact the same is true in arbitrary dimension and under much weaker
assumptions than general position:
\begin{proposition}
\label{prop:exact}
For any point set  $P\subset\reals^d$, $\im M\subseteq\ker \Delta$. If
$P$ has $d$ affinely independent points whose spanned hyperplane 
contains no other point of $P$, then the reversed inclusion also holds.
\end{proposition}

\begin{proof}
That $\im M\subset \ker\Delta$ (i.e., $\Delta\circ M =0$) is one direction of
Lemma \ref{lem:stress4}. To prove $\ker \Delta\subset\im M$, let $f\in
\reals^{\binom n 2}$ be such that $\sum_{ij\in C} w_{ij}f_{ij}=0$ 
for every circuit $C$. We want to find a motion $a=(a_1,\dots,a_n)$ for which 
$M(a)=f$. Without loss of generality assume that $p_1,\dots,p_d$ are
the claimed affinely independent points.
Let us fix a motion
$(a_1,\dots,a_{d})$ for the first $d$ points
satisfying the $\binom d2$ equations of the system
$Ma=f$ concerning these points.
This motion exists and is unique up to an arbitrary translational
and rotational component, because the complete graph on any affinely
independent point set is minimally infinitesimally
rigid.

By assumption, % The claim implies that 
$\{p_1,\dots,p_d,p_{i}\}$ is an affine basis for every $i>d$,
and hence the complete graph on it is again minimally infinitesimally rigid.
Thus, there is a motion
$(a'_1,\dots,a'_{d},a_i)$ 
satisfying the $\binom {d+1} 2$ corresponding
equations of the system
$Ma=f$. Adding translations and rotations
we can assume $a'_j=a_j$, $j=1,\dots,d$. So we have constructed a motion
$(a_1,\dots,a_n)$ which restricts to the one chosen for the $d$ first points
and which satisfies all the equations for pairs of points that include one of
the first $d$.

It remains to show that the equations are also
satisfied for the edges $kl$ with $k,l>d$.
We follow the same idea as in the proof of
Lemma~\ref{lem:equivalent-new}. Our assumption on the point set implies
that the unique circuit $C$ contained in $\{p_1,\dots,p_d,p_k,p_l\}$ uses both
$p_k$ and $p_l$. To simplify notation, assume that this circuit is 
$\{p_1,\dots,p_i,p_k,p_l\}$. By Lemma~\ref{lem:stress4},
there is a feasible motion $(a'_1,\dots,a'_{i},a'_k,a'_l)$. By translations
and rotations we assume $a_j=a'_j$ for $j=1,\dots,i$.
Observe now the
value of $\langle v-u, p_{j_1}-p_{j_2} \rangle$, for any of the points in the
circuit and for any vectors $v$ and $u$, depends only on the projection of
$v$ and $u$ to the affine subspace spanned by $C$. Since the complete
graphs on $\{p_1,\dots,p_i\}$, on
$\{p_1,\dots,p_i,p_k\}$ and on $\{p_1,\dots,p_i,p_l\}$ are 
minimally infinitesimally
rigid {\em when motions are restricted to 
that subspace}, we
conclude that the projections of $a_k$ and $a_l$ to that affine subspace
coincide with the projections of $a'_k$ and $a'_l$. In particular, 
$\langle a_l-a_k, p_l-p_k \rangle
=\langle a_l'-a_k', p_l-p_k \rangle = f_{kl}.$
\end{proof}

Hence, Lemma \ref{lem:equivalent-new} and its corollary,
Lemma~\ref{lem:meta-ppt},
hold in this generalized setting, with one equation per circuit instead of 
one equation per 4-tuple. The weakened general position assumption for $d$
of the points holds for every planar point set, since,
by Sylvester's theorem, any finite set of points
in the plane, not all on a single line,
has a line passing through only two of the points.

In dimension 3, however, the same is not true, and 
%some condition is necessary
%for the statement to hold. 
actually there are point sets for which $\im M \ne \ker \Delta$.
Consider the case of six points, three of them in
one line and three in another, with the two lines being skew (not parallel and
not crossing). These two sets of three points are the only two circuits
in the point set. In particular, $\ker \Delta$ has at most
codimension 2 in $\reals^{15}$, i.e., it has dimension at least~13.
On the other hand,
$\im M$ has at most the dimension of the reduced space of motions, $18-6=12$.

\section{Final Comments}

\label{final}

We describe some
open questions and ideas for further research.
The main questions related to this work are how to extend the
constructions from dimensions $1$ and $2$ to $3$ and higher, and how to
treat subsets in special
position in 2d. The
expectation is that this would give a coherent definition for
``pseudo-triangulations'' in higher dimensions.
Some ideas in this direction have been mentioned in
Section~\ref{towards}.

Is our choice of $f_{ij}$'s in Section~\ref{polytope} essentially unique?
The set of valid choices for a fixed point set is open in
$\reals^{\binom n 2}$. But, what if we restrict our attention to choices for
which, as in Theorem~\ref{thm:valid-new}, each $f_{ij}$ depends
only on the points $p_i$ and $p_j$, and not the rest of the
configuration? Observe that if we want this, then Theorem~\ref{thm:R}
provides an infinite set of conditions on the infinite set of unknowns 
$\{f(p,q) : p,q\in \reals^2\}$. It
follows from Lemma~\ref{lem:equivalent-new} that adding to a valid
choice $(f_{ij})_{i,j\in\{1,\dots,n\}}$ any vector 
$(\delta_{ij})_{i,j\in\{1,\dots,n\}}$ in the image of the rigidity map we still
  get a valid choice. And, of course, we can also scale any valid choice by a
  positive constant. This gives a half-space of valid choices of dimension 
${\binom n 2} +1$. Is this all of it?

% ------------------------
% SHOULD WE KEEP THE FOLLOWING PARAGRAPH? Personally, I don't think that the
% choices of parameters are so mysterious anymore. And perhaps the new, 
% transparent, proof
% works equally well for the wheels. I have not checked that.

It would also be interesting to see if there is a deeper reason
behind % the identity in
Lemma~\ref{lem:valid-2d}. % Theorem~\ref{thm:valid-new}. 
%For example, a more transparent proof of 
%the algebraic identities
%$R=2$ and $S=1$ in
%Lemmas~\ref{lem:R=S}--\ref{lem:invariant} would be desirable.
%These identities are in a sense trivial, as they can be proved
%by expanding all expressions and canceling terms (better
%with the aid of computer algebra software!).
We have actually been able to extend the identity
$\sum w_{ij} f_{ij}=1$
to a more general class of planar graphs than just the
complete graph on four vertices: to wheels (graphs of pyramids).
A wheel is a cycle with an additional vertex attached
to every vertex of the cycle.
For a wheel embedded
in the plane in general position,
 formulas that are quite
similar to~\eeqref{stress} in Lemma~\ref{lem:stresses}
define a self-stress $w_{ij}$
on its edges
which make the identity
$\sum w_{ij} f_{ij}=1$ true.
(Since the wheel is infinitesimally rigid 
and has $2n-2$ edges, the self-stress is unique up to a scalar factor.)
We take this as a hint that the identity % $R=1$ 
of Lemma~\ref{lem:valid-2d}
might be
an instance of a more general phenomenon which we don't fully understand.

\begin{figure}[htb]
  \centering
{\def\IPEfile{canonical.ipe}\input{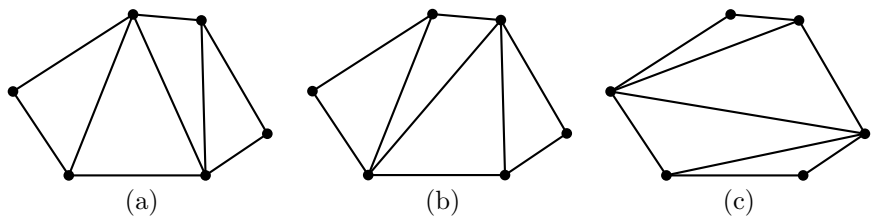}}
\caption{
(a)~The Delaunay triangulation of a point set in convex position.
(b-c)~The triangulations minimizing and maximizing
the objective function $(p_1,\dots,p_n)$ over the
ppt-polytope, respectively.
The triangulations in (b) and (c) are invariant under
affine transformations, whereas (a) is not.
}
\figlab{canonical}
\end{figure}

Another issue, with which we raised the introduction, is that having
a representation of combinatorial structures as vertices of
a polytope opens the way for selecting a particular
structure, by optimizing some linear functional over the polytope.
For example, the minimization of the objective function
with coefficient vector
$(|p_1|^2,\dots,|p_n|^2)$ 
%$(1,\dots,1)$ 
over the secondary polytope gives the
Delaunay triangulation. The opposite choice gives the furthest-site Delaunay
triangulation.
%What are the special ``canonical'' \emph{pointed pseudo-triangulations}
%that can be defined in this way?
The most natural choice of objective function for the polytope of pointed
pseudo-triangulations is $(p_1,\dots,p_n)$ or its opposite,
i.e., minimize
or maximize $\sum_{i} \langle p_i,v_i \rangle$ over all constrained expansions
which are tight on convex hull edges.
Even if, for points in convex position, our ppt-polytope is affinely
isomorphic to the secondary polytope, this functional on the ppt-polytope does
not, in general, give the Delaunay triangulation of those points,
see~Figure~\figref{canonical}.
%\complaint{\bf [I don't know which is max and which is min.,GR!!!]}
In fact, the result on the ppt-polytope is
invariant under affine transformations of the
point set, while the Delaunay triangulation is not.
% The investigation of 
The properties of the
pointed pseudo-triangulations that are defined in this way
await further studies.

\paragraph{Added in proof.} The second author, together with David Orden,
has extended the main construction of this paper to a simple polyhedron of
dimension $3n-3$ with a unique maximal bounded face whose vertices are 
\emph{all} the pseudo-triangulations of the point set.
Bounded edges correspond to either classical edge-flips or to the creation or destruction
of pointedness at a vertex by the deletion or inclusion of a single edge. The face poset
of this polyhedron is (essentially) the poset of all non-crossing graphs on the point set.

\endgroup

%\begin{figure}[htb]
%  \centering
%  \includegraphics{figure1.ps}
%  \caption{Including an [Encapsulated] Postscript file}
%  \label{no-ext}
%\end{figure}

% AFFILIATIONS
\section*{About the Authors}

G\"unter Rote is at the
Institut f\"ur Informatik,
Freie Universit\"at Berlin,
Taku\-stra{\ss}e 9, D-14195~Berlin,~Germany,
\textsl{rote@inf.fu-berlin.de}.
Francisco Santos is at the
Departamento de Matem\'aticas, Estad\'{\i}stica y Computaci\'on,
Universidad de Cantabria,
E-39005~Santander,~Spain,
\textsl{santos@matesco.unican.es}.
Ileana Streinu is at the Department of Computer Science, Smith College,
Northampton, MA~01063,~USA,
\textsl{streinu@cs.smith.edu}.

% AFFILIATIONS AND ACKS
\section*{Acknowledgments}

We thank Ciprian Borcea for pointing  out some relevant references in
the mathematical literature.
This work started at the Workshop on Pseudo-triangulations held at the McGill
University Bellairs Institute in Barbados, January 2001, partially funded by NSF grant
CCR-0104370, and continued during a visit of the third author to
the Graduiertenkolleg \emph{Combinatorics, Probability
and Computing} at Freie
Universit\"at Berlin,
supported by Deutsche Forschungsgemeinschaft, grant GRK~588/1.
Work by Ileana Streinu was partially supported by NSF RUI
Grants CCR-9731804 and CCR-0105507.
Work by Francisco Santos was partially supported
by grants PB97--0358 and BMF2001-1153
of the Spanish  Direcci\'on General de Ense\~nanza
Superior.

\end{document}